\documentclass[reqno]{amsart}
\usepackage[utf8x]{inputenc}
\usepackage{graphicx}
\usepackage{color}
\usepackage{float}
\usepackage{amsthm}
\usepackage{amssymb}
\usepackage{amsmath}
\usepackage{verbatim}
\usepackage{url}
\usepackage{tikz}
\usepackage{stmaryrd}

\usepackage{qtree}
\usepackage{bm}
\usepackage[a4paper,margin=1in,footskip=0.25in]{geometry}
\numberwithin{equation}{section}

\usetikzlibrary{chains,arrows,positioning}
\usetikzlibrary{matrix,arrows}

\theoremstyle{plain}
\newtheorem{theorem}{Theorem}[section]
\newtheorem{prop}[theorem]{Proposition}
\newtheorem{cor}[theorem]{Corollary}
\newtheorem{lemma}[theorem]{Lemma}
\newtheorem{conj}[theorem]{Conjecture}

\theoremstyle{definition}
\newtheorem{defn}[theorem]{Definition}
\newtheorem{example}[theorem]{Example}

\theoremstyle{remark}
\newtheorem{remark}[theorem]{Remark}

\title[Tropical critical points of the superpotential of a flag variety]{Tropical critical points of the superpotential of a flag variety}
\author{Jamie Judd}

\begin{document}

\begin{abstract}
In this paper we investigate the idea of a tropical critical point of the superpotential for the full flag variety of type $A$. Recall that associated to an irreducible representation of $G=SL_{n}(\mathbb{C})$ are various polytopes whose integral points parameterize a basis for the representation, e.g. the Gelfand-Zetlin polytope. Such polytopes can be constructed via the theory of geometric crystals by tropicalising a certain function, and in fact, the function involved coincides with the superpotential from the Landau-Ginzburg model for $G/B$ coming from mirror symmetry. In mirror symmetry a special role is played by the critical points of the superpotential, and motivated by this, we give a definition of the tropical critical points and use it to find a canonical point in each polytope. We then characterise the highest weights for which this tropical critical point is integral and therefore corresponds to a basis vector of the corresponding representation. Finally we give an interpretation of the tropical critical point by constructing a special vector in the representation using Borel-Weil theory and conjecturing a correspondence between this vector and the tropical critical point.
\end{abstract}
\maketitle
\setcounter{tocdepth}{1}
\tableofcontents
\section{Introduction}

Let $G$ be a simple complex algebraic group, $T$ a choice of maximal torus and $B \supset T$ a choice of Borel
subgroup, with opposite Borel subgroup $B_{-}$. Given $\lambda\in P^+$, a dominant integral weight of $T$, let $V_{\lambda}$ be the irreducible representation of $G$ with highest weight $\lambda$ and $V_{\lambda}^*$ the dual representation. Consider the case where
$\lambda=2 \rho$, the sum of the positive roots of G. There is a special vector in the representation $V_{\lambda}^*$ which has geometric origin and is defined as follows. First recall that Borel-Weil theory gives a geometric construction of $V_{\lambda}^*$ as $H^{0}(G/B,\mathcal{L}_{\lambda})$, where $\mathcal{L}_{\lambda}$ is the line bundle $$G \times^{B} \mathbb{C}_{-\lambda}=\{(g,x)\}/(g,x)\sim (gb,\lambda(b)x) \quad\text{(see [Sp]).}$$
If $\lambda=2 \rho$ then $\mathcal{L}_{2\rho}$  happens to be the anti-canonical bundle of $G/B$, so $V_{\lambda}^*$ is given by the global
sections of the anti-canonical bundle. Now, there exists a special meromorphic volume form $\omega$ on $G/B$, defined uniquely up to sign. This form was first introduced in [R] where it was defined as a
natural generalisation of the unique torus-invariant volume form on a torus inside a toric variety. It is
the meromorphic differential form on $G/B$ with simple poles exactly along the divisor given by the union
of all the Schubert divisors and all the opposite Schubert divisors, see [Lam, Section 2]. Similar volume
forms also appear more recently in work on mirror symmetry and cluster varieties, see [GHK, BMRS].
Now if we take the inverse of $\omega$, we get a special global section of the anti-canonical bundle of $G/B$, and
thus a distinguished vector in the representation $V_{2\rho}^*$. We would like to give an interpretation of this special vector.
\begin{example}
In the case $G=SL_{2}(\mathbb{C})$, we have $G/B\simeq \mathbb{P}^{1}$ via $\begin{pmatrix}g_{1}&g_2\\g_{3}&g_{4}\end{pmatrix}\mapsto x=\frac{g_1}{g_3}$. The Schubert divisors are given by $x=0$ and $x=\infty$ and the volume form $\omega$ is given by $\frac{dx}{x}$ which has a simple pole at $x=0$ and $x=\infty$. The representation $V_{2 \rho}^{*}$ is given by $\langle x^{2}\frac{\partial}{\partial x},x\frac{\partial}{\partial x},\frac{\partial}{\partial x}\rangle_{\mathbb{C}}$ and $\omega^{-1}$ is given by $x\frac{\partial}{\partial x}\in V_{2 \rho}^{*} $.\\
\end{example}

Restrict now to the case $G=SL_{n}(\mathbb{C})$. We will interpret this special section $\omega^{-1}$ using the mirror dual Landau-Ginzburg model for $G/B$. A Landau-Ginzburg model for the full flag variety of type $A$ was first introduced by Givental [G], in the form of
a regular function on a torus. Later Rietsch [R], inspired by [G] and work of Peterson [Pet97], gave a
construction of a Landau-Ginzburg model for partial flag varieties $G/P$ of general type, which in type $A$
is an extension of Givental’s function to a partial compactification of his torus. The Landau-Ginzburg model defined by Rietsch is a certain geometric object associated with the Langlands dual group $G^{\vee}$. It consists of a affine variety $Z\subset G^{\vee}$ and a holomorphic function $\mathcal{W}:Z \rightarrow  \mathbb{C}^{*}$. In the case of $G/B$, the variety $Z$ is given by the intersection of the Borel subgroup $B^{\vee}_{-}$ of $G^{\vee}$ with the big Bruhat double coset $B^{\vee}w_0B^{\vee}$ and the holomorphic function $\mathcal{W}:Z \rightarrow  \mathbb{C}^{*}$ is given by $$ u_1\bar{w}_0q^{-1}u_2\mapsto \chi(u_1)+\chi(u_2)\quad\text{where}\quad u_1,u_2\in U^{\vee},q\in T^{\vee}$$ Here $\chi$ is the character which takes the sum of the above diagonal entries and $\bar{w}_0$ is a particular representative in $G^{\vee}$ of the longest element of the Weyl group \footnote{In fact Rietsch used the decomposition $u_1q\bar{w}_0u_2 $, the decomposition we use here appears in [Ch]}. We also have that $Z$ fibers over the Langlands dual maximal torus $T^{\vee}$ via the map $$\textrm{hw}:Z\rightarrow T^{\vee}\quad\text{  given by  }\quad u_1\bar{w}_0q^{-1}u_2\mapsto q$$ For $q\in T^{\vee}$ the fiber $Z_{q}$ is isomorphic to a particular open subvariety of $G^{\vee}/B^{\vee}$ (in fact this open subvariety is the complement of the kind of divisor discussed above). Let the restriction of $\mathcal{W}$ to this fiber be $$\mathcal{W}_{q}:Z_{q} \rightarrow  \mathbb{C}^{*}.$$ This holomorphic function $\mathcal{W}$ first appeared in the work of [BK] on geometric crystals.
\begin{example}
In the case $G=SL_{2}(\mathbb{C})$, we have $Z=\{b=\begin{pmatrix}b_{1}&0\\b_{3}&b_{2}\end{pmatrix}\in B^{\vee}_{-}:b_3 \neq0 \}$ inside $G^{\vee}=PSL_{2}(\mathbb{C})$. This is isomorphic to $ T^{\vee}\times \mathbb{C}^*$ via $$ (\begin{pmatrix}q&\\&1\end{pmatrix},z)\mapsto\begin{pmatrix}1&z\\&1\end{pmatrix}\bar{w}_{0}\begin{pmatrix}q&\\&1\end{pmatrix}^{-1}\begin{pmatrix}1&q/z\\&1\end{pmatrix} =\begin{pmatrix}z/q&0\\1/q&1/z\end{pmatrix}$$ (We can think of $q\in T^{\vee}$ since $T^{\vee}\simeq \mathbb{C}^*$.) Hence each fiber $Z_{q}\simeq\mathbb{C}^{*} \subset \mathbb{P}^{1}\simeq G^{\vee}/B^{\vee}$, and we can see the function $\mathcal{W}_{q}:Z_{q} \rightarrow  \mathbb{C}^{*}$ is given by $z \mapsto z +\frac{q}{z}$.\\
\end{example}

Recall that associated to an irreducible representation of $G$ are various polytopes whose integral points parameterize a basis for the representation, e.g. the Gelfand-Zetlin polytope [GT]. In [BK], Berenstein and Kazhdan constructed such polytopes for each irreducible representation by applying tropicalisation to the function $\mathcal{W}:Z \rightarrow  \mathbb{C}^{*}$. In order to apply tropicalisation, we consider the variety $Z$ over the field of Puiseax series $\mathcal{K}$, along with the functions $\mathcal{W}:Z(\mathcal{K}) \rightarrow \mathcal{K}$ and $\mathrm{hw}:Z(\mathcal{K}) \rightarrow T^{\vee}(\mathcal{K})$. The field of Puiseax series is defined to be $$\mathcal{K}=\bigcup_{n=1}^{\infty}\mathbb{C}((t^{1/n}))$$ This consists of series $\sum_{p\geq q}c_{p}t^{p/M}$ where $q\in \mathbb{Z}$, $c_{p}\in \mathbb{C}$, $M\in \mathbb{Z}_{>0}$. Given an element $x(t)\in\mathcal{K}^*$ let $\mathrm{val}(x(t))\in \mathbb{Q}$ be the exponent appearing in the first non-zero term. Let $$\mathcal{K}_{>0}:=\{x(t)\in\mathcal{K}:\text{the coefficient of the lowest term of }x(t)\text{ is }\in \mathbb{R}_{>0}\}$$

There is a a well-defined notion, due to Lusztig [L], of the totally positive part of $Z(\mathcal{K})$, analogous to the case when the field is $\mathbb{C}$. Denote it by $Z(\mathcal{K}_{>0})$. If $\lambda\in P^+$ then $\lambda$ is a cocharacter of $T^{\vee}$ and we can define $t^{\lambda}\in T^{\vee}(\mathcal{K})$. Define $$Z_{t^\lambda}(\mathcal{K}_{>0}):=\{z\in Z(\mathcal{K}_{>0}):\mathrm{hw}(z)=t^\lambda\}$$ and $$Z_{t^\lambda}(\mathcal{K}_{>0})^{+} := \{z\in Z_{t^\lambda}(\mathcal{K}_{>0}):\mathrm{val}(\mathcal{W}(z))\geq0\}$$

Following Berenstein and Kazhdan, for each reduced expression $\textbf{i}$ of $w_0$ we will define a valuation map (see equation 4.1) $$\nu^{\vee}_{\lambda,\textbf{i}}:Z_{t^\lambda}(\mathcal{K}_{>0})^+\longrightarrow \mathbb{R}^{\resizebox{0.06\hsize}{!}{$dim(G/B)$}}$$ which has the property that the closure of the image of $\nu^{\vee}_{\lambda,\textbf{i}}$ is the string polytope associated to $\textbf{i}$. This is a polytope introduced by Littlemann [Lit], whose integral points parameterise a basis for $V_{\lambda}$.

\begin{example}
Take $G=SL_{2}(\mathbb{C})$. Let $\lambda\in \mathbb{Z}_{>0}$ and think of it as the dominant integral weight $\lambda \alpha_1$. Then we get that $$Z_{t^\lambda}(\mathcal{K}_{>0})^{+} \simeq \{ z\in\mathcal{K}_{>0}:\textrm{val}(z+\frac{t^{\lambda}}{z})\geq0\} = \{ z=at^c+...\in\mathcal{K}_{>0}:\min(c,\lambda-c)\geq0\}$$ In this case, it turns out the map $\nu^{\vee}_{\lambda,\textbf{i}}$ will be just given by the valuation $\mathrm{val}:\mathcal{K}_{>0}\rightarrow\mathbb{Q}$ (for any $\textbf{i}$). Hence the image is $\{c\in \mathbb{Q}:\min(c,\lambda-c)\geq0\}$ which has closure the 1-dimensional polytope $[0,\lambda]$. The integral points are $\{0,1,2,...,\lambda-1,\lambda\}$.\\
\end{example}

There is also another way to construct such polytopes using the theory of Newton-Okounkov bodies. This construction takes a line bundle $\mathcal{L}$ on a projective variety $X$, a section $\tau\in H^0(X,\mathcal{L})$ and a valuation $val:\mathbb{C}(X)\rightarrow \mathbb{Z}^{dim(X)}$ and produces a convex body which in the simplest cases is just the convex hull of the image of the map $H^0(X,\mathcal{L})\rightarrow\mathbb{Z}^{dim(X)}\subset \mathbb{R}^{dim(X)}$ sending $\sigma\mapsto val(\frac{\sigma}{\tau})$. Given a choice of $\lambda\in P^+$ and reduced expression $\textbf{i}$ for $w_0$, Kaveh [Ka] defined a choice of valuation map $$\nu_{\lambda,\textbf{i}}:H^{0}(G/B,\mathcal{L}_{\lambda})\rightarrow \mathbb{R}^{\resizebox{0.06\hsize}{!}{$dim(G/B)$}}$$ whos associated Newton-Okounkov body is exactly the string polytope for the reduced expression $\textbf{i}$.

\begin{example}
Continuing the example $G=SL_{2}(\mathbb{C})$, we have $G/B=\mathbb{P}^{1}$ and $\mathcal{L}_{\lambda}=\mathcal{O}(\lambda)$. Then $$V_{\lambda}^{*}=H^{0}(\mathbb{P}^{1},\mathcal{O}(\lambda))\simeq \langle x^{\lambda},...,x^{1},1\rangle_{\mathbb{C}}$$ In this case $\nu_{\lambda,\textbf{i}}$ is just the natural valuation on $\mathbb{C}(x)$, so we get that the image is $\{0,1,2,...,\lambda\}$ with convex hull $[0,\lambda]$.\\
\end{example}

Now we return our attention to giving an interpretation of the distinguished vector discussed earlier. This will be done by considering the critical points of the function $\mathcal{W}$. In the context of the Landau-Ginzburg model for a partial flag variety $G/P$, the critical points have an important interpretation as Rietsch showed that the fiber-wise critical points of $\mathcal{W}$ describe the quantum cohomology $QH^{*}(G/P)$ of $G/P$.
\begin{example}
In our running example, the fiber-wise critical points of $\mathcal{W}$ are given by $\frac{\partial{\mathcal{W}}}{\partial z}=1-\frac{q}{z^{2}}=0$, and the quantum cohomology ring of $G/B=\mathbb{P}^{1}$ is given by $\mathbb{C}[z,q]/(z^{2}-q)$, where $z$ represents the hyperplane cohomology class.\\
\end{example}
Motivated by the fact that the critical points of $\mathcal{W}_{q}:Z_{q}\rightarrow\mathbb{C} $ have an important interpretation in terms of the quantum cohomology of the original flag variety, we study a tropicalised version of the critical points by considering the critical points of the function $\mathcal{W}_{t^{\lambda}}:Z_{t^{\lambda}}(\mathcal{K}) \rightarrow \mathcal{K}$. In [R], Rietsch showed that $\mathcal{W}_{q}$ has a unique critical point in the totally positive part of $Z_{q}$, and in Theorem \ref{tcp} we prove a similar result: that there exists a unique critical point of $\mathcal{W}_{t^{\lambda}}$ in the totally positive part $Z_{t^{\lambda}}(\mathcal{K}_{>0})$. We denote this point by $p_{\lambda}$ and also show that in fact $p_{\lambda}\in Z_{t^{\lambda}}(\mathcal{K}_{>0})^+$.

We are especially interested in the case when $p_{\lambda}$ lies in $Z(\mathbb{C}((t)))$, in which case we say $p_{\lambda}$ is integral. In section 6, we characterise the $p_{\lambda}$ which are integral.

\begin{example}
In the example $G=SL_{2}(\mathbb{C})$, the critical points are given by $\{ z=at^c+...\in\mathcal{K}_{>0}:1-\frac{t^{\lambda}}{z^{2}}=0\}$, i.e. $z=t^{\frac{\lambda}{2}}$. This point will be integral whenever $\lambda$ is even.\\
\end{example}

Thus, given $\textbf{i}$ a reduced expression of $w_0$, we have two maps
\begin{align*}
\nu^{\vee}_{\lambda,\textbf{i}}:Z_{t^\lambda}(\mathcal{K}_{>0})^+\longrightarrow \mathbb{R}^{\resizebox{0.06\hsize}{!}{$dim(G/B)$}}\\
\nu_{\lambda,\textbf{i}}:H^{0}(G/B,\mathcal{L}_{\lambda})\longrightarrow \mathbb{R}^{\resizebox{0.06\hsize}{!}{$dim(G/B)$}}
\end{align*}
both of whose images lie in the string polytope for $\textbf{i}$.
Our interpretation of the special section $\omega^{-1}\in H^{0}(G/B,\mathcal{L}_{2\rho})$ is via the following conjecture:
\begin{conj}
The image $\nu_{2\rho,\textbf{i}}(\omega^{-1})=\nu^{\vee}_{2\rho,\textbf{i}}(p_{2\rho})$ for all reduced expressions $\textbf{i}$ of $w_0$.
\end{conj}

\begin{example}
In the example, $\nu^{\vee}_{\lambda,\textbf{i}}(p_{\lambda})$ is given by the point $\lambda/2 \in [0,\lambda]$. When $\lambda=2\rho$ the polytope is $[0,2]$ and the valuation of the critical point gives the point $1\in[0,2]$. Now consider the Newton-Okounkov construction in the case $\lambda=2 \rho$. The image of $\omega^{-1}=x\frac{\partial}{\partial x}$ under $\nu_{\lambda,\textbf{i}}$ gives the point $1\in[0,2]$, which is in agreement.
\end{example}

We also generalise this as follows. For each $\lambda$ such that $p_{\lambda}$ is integral, we will define a special section $\omega_{\lambda}^{-1}\in H^{0}(G/B,\mathcal{L}_{\lambda})$ which generalises the section $\omega^{-1}\in H^{0}(G/B,\mathcal{L}_{2\rho})$ defined above. Then we make a similar conjecture:

\begin{conj}
Given $\lambda$ such that $p_{\lambda}$ is integral, we have $\nu_{\lambda,\textbf{i}}(\omega_{\lambda}^{-1})=\nu^{\vee}_{\lambda,\textbf{i}}(p_{\lambda})$ for all reduced expressions $\textbf{i}$ of $w_0$.
\end{conj}

\textbf{Acknowledgments}: I would especially like to thank my supervisor Konni Rietsch for suggesting this idea to me and for helpfully guiding me through it. I would like to acknowledge Victor Ginzburg who originally posed the problem of interpreting the element of $V_{2 \rho}^{*}$ defined by the special global section of the anticanonical bundle. I would also like to acknowledge Nicholas Shepherd-Barron whose discussion with Konni led to the idea of using the ideal fillings.\\

\section{Background and Notation}

Let $G=SL_{n}(\mathbb{C})$. Let $B$ and $B_{-}$ be the upper triangular and lower triangular matrices in $G$, let $U$ and $U_{-}$ be the strictly upper triangular and strictly lower triangular matrices in $G$, and let $T$ be the diagonal matrices in $G$. Let $G^{\vee}=PGL_{n}(\mathbb{C})$ and let $B^{\vee},B_{-}^{\vee},U^{\vee},U_{-}^{\vee}$ and $T^{\vee}$ be the corresponding subgroups in $G^{\vee}$. We have $T \subset T_{GL_{n}}$ and $T_{GL_{n}^{\vee}} \twoheadrightarrow T^{\vee}$ where $T_{GL_{n}}$ and $T_{GL_{n}^{\vee}}$ are the diagonal matrices in $GL_{n}$ and its Langlands dual group $GL_{n}^{\vee}$. For $i=1,...,n$ let $\epsilon_{i}$ and $\mu_{i}$ be the standard characters of $T_{GL_{n}}$ and $T_{GL_{n}^{\vee}}$ respectively. Then we have $$X^{*}(T)=X_{*}(T^{\vee})=\langle \epsilon_{1},...,\epsilon_{n} \rangle / (\sum \epsilon_{i})$$ and $$X_{*}(T)=X^{*}(T^{\vee})= \{\sum c_{i}\mu_{i} : \sum c_{i} =0\}$$
Let $\alpha_{ij}=\epsilon_{i}-\epsilon_{j}\in X^{*}(T)$ and $\alpha_{ij}^{\vee}=\mu_{i}-\mu_{j}\in X^{*}(T^{\vee})$. The roots of $G$ are $\Delta_{G}=\{\alpha_{ij}:i \neq j\}$ and the positive roots are $\Delta_{G}^{+}=\{\alpha_{ij}:i < j\}$. Let $P^+ :=\{\lambda \in X^{*}(T):\langle\lambda,\alpha^{\vee}_{ij}\rangle\geq 0\quad \forall i<j\} $ be the set of dominant integral weights. Let $\rho$ denote half the sum of the positive roots of $G$. Let $I=\{1,2,...,n-1\}$. Given $i \in I$ we let $\alpha_{i}:=\alpha_{ii+1}$, and thus we have that $\{\alpha_{i}:i \in I\}$ are the simple roots of $G$ (and the simple coroots of $G^{\vee}$). We also have the fundamental weights of $G$ given by $\omega_{i}:=\epsilon_{1}+...+\epsilon_{i}$ and recall the equality $\rho=\sum \omega_{i}$. Finally, let $W=S_{n}$ be the Weyl group and let $\{s_{i}=(i,i+1):i \in I\} $ be the simple reflections.

For $i \in I$, let $\phi_{i}$ be the homomorphism $SL_{2} \rightarrow G$ corresponding to the $i^{th}$ simple root of $G$. Given $i \in I$ define $$\textbf{x}_{i}:\mathbb{G}_{a} \rightarrow G \quad \text{  to be  }\quad z \mapsto \phi_{i}\begin{pmatrix}1&z\\0&1\end{pmatrix}$$  $$\textbf{y}_{i}:\mathbb{G}_{a} \rightarrow G \quad \text{  to be  }\quad z \mapsto \phi_{i}\begin{pmatrix}1&0\\z&1\end{pmatrix}$$
 $$\textbf{x}_{-i}:\mathbb{G}_{a} \rightarrow G \quad \text{  to be  }\quad z \mapsto \phi_{i}\begin{pmatrix}z^{-1}&0\\1&z\end{pmatrix}$$
and note that $\textbf{x}_{-i}(z)=\textbf{y}_{i}(z)\alpha_{i}^{\vee}(z^{-1})$. For a simple reflection $s_{i} \in W$, define $$\dot{s}_{i}=\textbf{x}_{i}(-1)\textbf{y}_{i}(1)\textbf{x}_{i}(-1)\in G$$ For general $w \in W$, a representative $\dot{w} \in G$ is defined by $\dot{w}=\dot{s}_{i_{1}}\dot{s}_{i_{2}}...\dot{s}_{i_{m}}$, where $s_{i_{1}}s_{i_{2}}...s_{i_{m}}$ is a (any) reduced expression for $w$. The length of a reduced expression for $w$ is denoted by $l(w)$. Let $N=l(w_{0})$ where $w_{0}$ is the longest element of $W$.

Similarly for $G^{\vee}$ we have a homomorphism $\phi_{i}^{\vee}:PGL_{2} \rightarrow G^{\vee}$ for each $i \in I$ and in the same way we can define maps $\textbf{x}_{i}^{\vee}, \textbf{y}_{i}^{\vee}$, $\textbf{x}_{-i}^{\vee}$ from $\mathbb{G}_{a}$ to $G^{\vee}$. Also, in the same way as above, we can define a representative of $w$ in $G^{\vee}$ which we denote by $\bar{w}$.

For $P \supset B$ a parabolic subgroup of $G$, let $I_{P}=\{i \in I \vert \dot{s}_{i} \in P\}$ and $I^{P}$ be its complement in $I$. Note that if we write $I^{P}=\{n_{1},...,n_{k}\}$ where $n_{1} < n_{2}<...<n_{k}$ then $G/P$ is isomorphic to the partial flag variety $Fl(n_{1},n_{2},..,n_{k},n)$. Let $$W_{P}:=\langle s_{i} : i \in I_{P}\rangle$$ and let $w_{P}$ be the longest element in the subgroup $W_{P}$. Let $N_{P}=N-l(w_{P})$.

Given a parabolic subgroup $P$ with Levi subgroup of $L$, let $\Delta_{+}^{L}$ (respectively $\Delta_{-}^{L}$) the roots of $L$ which are positive (respectively negative), and define $$\lambda_{P}:=\sum\limits_{\Delta_{+}^{G}\backslash \Delta_{+}^{L}} \alpha$$ For example $\lambda_{B}=2 \rho$.

We also define the open Richardson varieties as follows. For $v,w\in W$ with $v \leq w$ (Bruhat order) define $$\mathcal{R}_{v,w}:=(B_{-} \dot{v} B \cap B \dot{w} B)/B \subset G/B$$ It is known that $\mathcal{R}_{v,w}$ is smooth and irreducible of dimension $l(w)-l(v)$. Also let $$\mathcal{R}^{\vee}_{v,w}:=B^{\vee} \backslash (B^{\vee} \bar{v} B^{\vee}_{-} \cap B^{\vee} \bar{w} B^{\vee}) \subset B^{\vee}\backslash G^{\vee}$$

We now introduce a particular subvariety $Z \subset G^{\vee}$ which appeared in [R] as a Landau-Ginzburg model and in [BK] as a decorated geometric crystal. Let $$Z:=B_{-}^{\vee} \cap B^{\vee} \bar{w}_{0}B^{\vee}$$
We equip $Z$ with three maps, the ``highest weight map'' $\mathrm{hw}:Z \rightarrow T^{\vee}$, the ``superpotential'' $\mathcal{W}:Z \rightarrow\mathbb{C}^{*}$ and the ``weight map'' $ \mathrm{wt}:Z \rightarrow T^{\vee}$. Firstly we define the highest weight map. Any element of $Z$ can be written uniquely as $u_{1}\bar{w}_{0}q^{-1}  u_{2}$ with $u_1,u_2\in U^{\vee}$ and $q\in T^{\vee}$, so we can define $$\mathrm{hw} :Z \rightarrow T^{\vee} \quad \text{ sending }\quad u_{1}\bar{w}_{0}q^{-1}  u_{2} \mapsto q$$ Let $Z_{q}$ be the fiber of $\mathrm{hw}$ over $q\in T^{\vee}$.
\begin{remark}
The map $Z \rightarrow B^{\vee} \backslash G^{\vee}$ which sends $b \mapsto B^{\vee} b$ gives an isomorphism $Z_{q}\simeq \mathcal{R}^{\vee}_{e,w_{0}}$. These combine to give an isomorphism $Z \simeq T^{\vee} \times\mathcal{R}^{\vee}_{e,w_{0}}$.
\end{remark}
\noindent Now define $\chi:U^{\vee}\rightarrow \mathbb{C}$ by $$\chi(u)=\sum\limits_{i=1}^{n-1}u_{i\,i+1}$$
and define the superpotential to be the map $$\mathcal{W}:Z \rightarrow\mathbb{C}^{*}\quad\text{  sending  }\quad u_{1}\bar{w}_{0}q^{-1}  u_{2} \mapsto \chi (u_{1})+\chi (u_{2})$$
The weight map is defined to be the map $$ \mathrm{wt}:Z \rightarrow T^{\vee}\quad\text{  sending  }\quad b \mapsto \bar{w}_0 pr(b)^{-1}\bar{w}_{0}^{-1}$$ where $pr$ is the projection $B_{-}^{\vee}\rightarrow B_{-}^{\vee}/U^{\vee}_{-}=T^{\vee}$.

\begin{example}
For $G=SL_3$ these maps are: $$\mathrm{hw}:\begin{pmatrix}b_1&&\\b_2&b_3&\\b_4&b_5&b_6\end{pmatrix}\mapsto\begin{pmatrix}\frac{1}{b_4}&&\\&\frac{b_4}{b_2b_5-b_3b_4}&\\&&\frac{b_2b_5-b_3b_4}{b_1b_3b_6}\end{pmatrix} $$
$$\mathcal{W}:\begin{pmatrix}b_1&&\\b_2&b_3&\\b_4&b_5&b_6\end{pmatrix}\mapsto \frac{b_2+b_5}{b_4}+\frac{b_1b_5+b_2b_6}{b_2b_5-b_3b_4}$$
$$\mathrm{wt}:\begin{pmatrix}b_1&&\\b_2&b_3&\\b_4&b_5&b_6\end{pmatrix}\mapsto\begin{pmatrix}\frac{1}{b_6}&&\\&\frac{1}{b_3}&\\&&\frac{1}{b_1}\end{pmatrix} $$
\end{example}

\begin{remark}
It is shown in [R] that the (fiber-wise) critical points of the superpotential trace out the open part of the Peterson variety, whose coordinate ring is isomorphic to the quantum cohomology ring of the Langlands dual flag variety $G/B$. See [R] for more details.\\
\end{remark}

Finally we introduce the twist map which will be needed later. Berenstein and Zelevinsky define an isomorphism called the twist map $$\eta^{w_{0},e}:B^{\vee}_{-}\cap U^{\vee}\bar{w}_{0}U^{\vee} \longrightarrow U^{\vee}\cap B^{\vee}_{-}\bar{w}_{0}B^{\vee}_{-}\quad\text{  sending  }\quad b \mapsto [(\bar{w}_{0}b^T)^{-1}]_{+} $$ where $b^T$ is the transpose of $b$ and $[g]_{+}$ is defined by the decomposition $g=[g]_{-}[g]_{0}[g]_{+}$ with $[g]_{-}\in U^{\vee}_{-}$, $ [g]_{0}\in T^{\vee}$, $[g]_{+}\in U^{\vee}$.

Note that we also have an isomorphism $$\Phi:T^{\vee}\times(  U^{\vee}\cap B^{\vee}_{-}\bar{w}_{0}B^{\vee}_{-}) \simeq Z\quad\text{  via  }\quad(q,u)\mapsto u\bar{w}_{0}q^{-1} u_2$$ where $u_2$ is the unique element in $ U^{\vee}$ making this lie in $Z$.
\section{Toric charts}
We now introduce the definition of a positive variety, and then we will give the variety $Z$ and the partial varieties $G/P$ the structure of positive varieties.

\begin{defn}\cite{BK}
Given an algebraic torus $S$ we say a rational function on $S$ is \textit{positive} if, in a coordinate system given by a set of characters of $S$, it can be written as a ratio of two polynomials with positive integral coefficients. A rational function $f:S \rightarrow S'$ between two algebraic tori is positive if $\chi \circ f$ is positive for all $\chi \in X^{*}(S')$. A \textit{toric chart} on an algebraic variety $Y$ is a birational isomorphism $\theta:S \rightarrow Y$. Two toric charts $\theta:S \rightarrow Y$ and $\theta':S' \rightarrow Y$ are \textit{positively equivalent} if $(\theta)^{-1} \circ \theta'$ and $(\theta')^{-1} \circ \theta$ are both positive. An equivalence class $\Theta_{Y}$ of toric charts is called a \textit{positive structure} on $Y$, and the pair $(Y,\Theta_{Y})$ is called a positive variety.
\end{defn}

First we define some toric charts on $Z$ indexed by reduced expressions for $w_{0}$. Given a reduced expression $\textbf{i}=(i_{1},...,i_{N})$ for $w_{0}$,
define $$\textbf{x}^{\vee}_{-\textbf{i}}:(\mathbb{C}^{*})^{N} \rightarrow B^{\vee}_{-}\cap U^{\vee}\bar{w}_{0}U^{\vee}  \quad\text{ by }\quad (z_{1},..,z_{N}) \mapsto \textbf{x}_{-i_{1}}^{\vee}(z_{1})...\textbf{x}_{-i_{N}}^{\vee}(z_{N})$$ Then define $$\tilde{\textbf{x}}_{-\textbf{i}}^{\vee}:T^{\vee}\times (\mathbb{C}^{*})^{N} \longrightarrow Z\quad\text{ by }\quad (q,z) \mapsto \Phi(q,\eta^{w_{0},e}(\textbf{x}^{\vee}_{-\textbf{i}}(z)))$$ The toric charts $\tilde{\textbf{x}}_{-\textbf{i}}^{\vee}$ on $Z$ are all positively equivalent, so they give a positive structure on $Z$, which we denote $\Theta_{Z}$. This follows from [BZ, Prop 4.5, Thrm 4.7]. See also [BK, 3.26].

Also observe that for a fixed $q \in T^{\vee}$ the restriction of any of these charts to $\{q\}\times (\mathbb{C}^{*})^{N} $ gives a toric chart on $Z_{q}$. Hence in any of these charts $\mathrm{hw}$ is given by the natural projection $T^{\vee}\times (\mathbb{C}^{*})^{N}  \rightarrow T^{\vee}$.

\begin{remark}
The union of the charts $\tilde{\textbf{x}}_{-\textbf{i}}^{\vee}\vert_{\{q\}\times (\mathbb{C}^{*})^{N} }$ as $\textbf{i}$ varies over reduced expressions for $w_{0}$ covers $Z_{q}$ up to at least codimension two. This follows from [R1, Lemma 3.7].
\end{remark}

\begin{example}
Let $G=SL_3$ and $\textbf{i}=(212)$. We have $$\tilde{\textbf{x}}_{-\textbf{i}}^{\vee}(q,z)=\begin{pmatrix}1&z_3&z_2\\0&1&z_1+\frac{z_2}{z_3}\\0&0&1\end{pmatrix}.\bar{w_0}.\begin{pmatrix}q_1^{-1}&&\\&q_2^{-1}&\\&&q_3^{-1}\end{pmatrix}.\begin{pmatrix}1&\frac{q_1 z_3}{q_2z_2}&\frac{q_1}{q_3z_1z_3}\\0&1&\frac{q_2 (z_2+z_1z_3)}{q_3z_1z_3^2}\\0&0&1\end{pmatrix}$$
So we can compute $\mathcal{W}$ in this chart to be
   $$(q,z) \mapsto z_3+z_1+\frac{z_2}{z_3}+\frac{q_1 z_3}{q_2z_2}+\frac{q_2 (z_2+z_1z_3)}{q_3z_1z_3^2}$$
\end{example}

On the other side we will need certain toric charts on $G/P$. Firstly, if $\textbf{i}=(i_{1},...,i_{N})$ is a reduced expression for $w_{0}$,
define $$\textbf{y}_{\textbf{i}}:S \rightarrow G/B \quad \text{   by   } \quad(x_{1},...,x_{N}) \mapsto \textbf{y}_{i_{1}}(x_{1})...\textbf{y}_{i_{N}}(x_{N})B$$
To define the toric charts on $G/P$ we need the definition of a \textit{positive subexpression}. The reference for this is [De] (see also [R,§6]).
 Let $\textbf{i}=(i_{1},...,i_{N})$ be a reduced expression for $w_{0}$, then given $v\in W$ there is a unique subexpression $(i_{j_{1}},..,i_{j_{t}})$ such that: $$s_{i_{j_{1}}}...\,s_{i_{j_{t}}}=v$$ $$(s_{i_{j_{1}}}...\,s_{i_{j_{l}}})s_{i_{h}}>s_{i_{j_{1}}}...\,s_{i_{j_{l}}} \quad \text{for all }\quad j_{l}<h\leq j_{l+1} \quad \text{and }\quad 1\leq l\leq t$$
It is a reduced expression and called the \textit{positive subexpression} for $v$. The algorithm to compute it can be described as follows:
Let $u=v$. For $t=N,...,1$; if $us_{i_{t}}<u$ put a circle around $s_{i_{t}}$ and replace $u$ by $us_{i_{t}}$. Then the circled factors give the positive subexpression.
\begin{example}
Suppose $G=GL_{5},I_{P}=\{2,3\},\textbf{i}=(4321432434)$, then the positive subexpression for $w_{P}=(232)=(323)$ is given by $43214\textcircled{3}\textcircled{2}4\textcircled{3}4$.
\end{example}
\noindent Now let $P \supset B$ be a parabolic subgroup of $G$ and let $\textbf{i}=(i_{1},...,i_{N})$ a reduced expression for $w_{0}$. Let $\textbf{j}=(i_{j_{1}},..,i_{j_{l(w_{P})}})$ be the positive subexpression for $w_{P}$. Let $J=\{j_{1},..,j_{l(w_{P})}\}$ and let the complement of $J$ be $\{k_{1},..,k_{N_{P}}\}$. Define $$\textbf{y}^{P}_{\textbf{i}}:(\mathbb{C}^{*})^{N_{P}} \rightarrow G/P \quad\text{   by   }\quad (u_{k_{1}},...,u_{k_{N_{P}}}) \mapsto g_{1}...g_{N}P \quad \text{ where }\quad g_{t}=\begin{cases} \textbf{y}_{i_{t}}(u_{t}) \quad t \notin J \\ \dot{s}_{i_{t}} \quad\quad \,\,\,\, t \in J \end{cases}$$
These charts are all positively compatible so define a positive structure on $G/P$. This follows from [R,§7].\\

\section{Polytopes associated to $\lambda\in P^+$}
In this section we will associate a polytope to a choice of dominant integral weight and reduced expression of $w_0$, in two different ways.

First we define the totally positive part of a positive variety over the field of Puiseax series. The field of Puiseax series is defined to be $$\mathcal{K}=\bigcup_{n=1}^{\infty}\mathbb{C}((t^{1/n}))$$ This consists of series $\sum_{p\geq q}c_{p}t^{p/M}$ where $q\in \mathbb{Z}$, $c_{p}\in \mathbb{C}$, $M\in \mathbb{Z}_{>0}$. Given an element $x(t)\in\mathcal{K}^*$ let $\mathrm{val}(x(t))\in \mathbb{Q}$ be the exponent appearing in the first non-zero term. We define a map $\mathrm{val}:(\mathcal{K}^*)^m\rightarrow \mathbb{R}^m$ sending $(x_1,...,x_m)\mapsto (\mathrm{val}(x_1),...,\mathrm{val}(x_m))$. Let $$\mathcal{K}_{>0}:=\{x(t)\in\mathcal{K}:\text{the coefficient of the lowest term of }x(t)\text{ is }\in \mathbb{R}_{>0}\}$$

Given any variety $Y$ we can consider it over the field $\mathcal{K}$, which we denote by $Y(\mathcal{K})$. If $S=(\mathbb{C}^{*})^m$ is an algebraic torus we can consider $S(\mathcal{K})$ and define the totally positive part to be $S(\mathcal{K}_{>0}):=(\mathcal{K}_{>0})^m$. Now given a positive variety $(Y,\Theta_{Y})$, define the totally positive part of $Y(\mathcal{K})$ to be $Y(\mathcal{K}_{>0}):=\theta(S(\mathcal{K}_{>0}))$ for any $\theta:S\rightarrow Y$ in $\Theta_{Y}$. This is independant of the choice of chart since the transition maps are positive, so therefore preserve $S(\mathcal{K}_{>0})$.

Now let $\lambda\in P^+$. Since $\lambda$ is a cocharacter of $T^{\vee}$ we can define $t^{\lambda}\in T^{\vee}(\mathcal{K})$. Define $$Z_{t^\lambda}(\mathcal{K}_{>0}):=\{z\in Z(\mathcal{K}_{>0}):\mathrm{hw}(z)=t^\lambda\}$$ and $$Z_{t^\lambda}(\mathcal{K}_{>0})^{+} := \{z\in Z_{t^\lambda}(\mathcal{K}_{>0}):\mathrm{val}(\mathcal{W}(z))\geq0\}$$

Following the ideas of Berenstein and Kazhdan, we make the following definition. Given $\textbf{i}=(i_1,...,i_N)$ a reduced expression of $w_0$, define $$\nu^{\vee}_{\lambda,\textbf{i}}:Z_{t^\lambda}(\mathcal{K}_{>0})^+\longrightarrow \mathbb{R}^N$$
via the composition \begin{equation}Z_{t^\lambda}(\mathcal{K}_{>0})^+\xrightarrow{(\tilde{\textbf{x}}_{-\textbf{i}}^{\vee})^{-1}} T^{\vee}(\mathcal{K}_{>0})\times (\mathcal{K}_{>0})^N \xrightarrow{\mathrm{pr}_2}(\mathcal{K}_{>0})^N \xrightarrow{\mathrm{val}} \mathbb{R}^N\end{equation}

Then we have the following theorem
\begin{theorem}
The closure of the image of $\nu^{\vee}_{\lambda,\textbf{i}}$ is the string polytope $\mathrm{String}_{\bf{i}}(\lambda)$.
\end{theorem}

This follows from [Ch] and the description of the string polytope in [Lit]. To see this we first introduce tropicalisation. Given a positive morphism between two algebraic tori we define its tropicalisation  by replacing $+$ by $\mathrm{min}$ and replacing $\times$ by $+$.
\begin{example}
 Suppose $\phi:(\mathcal{K}^{*})^{4}\rightarrow\mathcal{K}^{*}$ is the map $\frac{x_{1}x_{2}^{2}+3x_{2}x_{3}^{5}}{x_{2}x_{4}}$ then $\phi^{t}:\mathbb{R}^{4}\rightarrow \mathbb{R}$ is given by $(a_{1},a_{2},a_{3},a_{4})\mapsto \mathrm{min}\{a_{1}+2a_{2},a_{2}+5a_{3}\}-\mathrm{min}\{a_{2}+a_{4}\}$.
\end{example}
Then consider $\mathcal{W}_{\lambda,\textbf{i}}:(\mathcal{K}^*)^N \rightarrow \mathcal{K}^*$ defined by $z \mapsto \mathcal{W}(\tilde{\textbf{x}}_{-\textbf{i}}^{\vee}(t^\lambda,z))$ and its tropicalsation $\mathcal{W}_{\lambda,\textbf{i}}^t:\mathbb{R}^N\rightarrow \mathbb{R}$. We can describe the image of $\nu^{\vee}_{\lambda,\textbf{i}}$ as $\{a\in\mathbb{Q}^N:\mathcal{W}_{\lambda,\textbf{i}}^t(a)\geq0\}$ since $\mathrm{val}\circ\mathcal{W}_{\lambda,\textbf{i}}=\mathcal{W}_{\lambda,\textbf{i}}^t\circ\mathrm{val}$. Hence the closure of the image of $\nu^{\vee}_{\lambda,\textbf{i}}$ is $\{a\in\mathbb{R}^N:\mathcal{W}_{\lambda,\textbf{i}}^t(a)\geq0\}$.

From the definition of $\tilde{\textbf{x}}_{-\textbf{i}}^{\vee}$, we have that $\mathcal{W}_{\lambda,\textbf{i}}$ is given by $z=(z_1,...,z_N) \mapsto \chi(\eta^{w_0,e}(\textbf{x}_{-\textbf{i}}^{\vee}(z)))+\chi(u_2)$ where $u_2\in U^{\vee}$ is the unique element such that $\eta^{w_0,e}(\textbf{x}_{-\textbf{i}}^{\vee}(z)).\bar{w}_0.t^{-\lambda}.u_2$ lies in $B_{-}$. By [Ch] Proposition 5.2.5 this is equal to $$\chi(\eta^{w_0,e}(\textbf{x}_{-\textbf{i}}^{\vee}(z)))+\sum_{k=1}^{N}\alpha_{i_k}^{\vee}(t^{\lambda})z_k^{-1}\prod_{j=k+1}^{N}z_{j}^{-\langle \alpha_{i_j},\alpha_{i_k}^{\vee}\rangle}$$ Then tropicalising this expression and setting it $\geq0$ gives the conditions cutting out the string polytope given in [Lit] Proposition 1.5. Also see [Ch] Proposition 7.2.1.\\

\begin{example}
Let $G=SL_3$ and $\textbf{i}=(212)$. From Example 3.3, we get that $\mathcal{W}_{\lambda,\textbf{i}}^{t}$ is given by $$(c_1,c_2,c_3)\mapsto \min\{c_3,c_1,c_2-c_3,\lambda_1-\lambda_2+c_3-c_2,\lambda_2-\lambda_3+c_2-c_1-2c_3,\lambda_2-\lambda_3-c_3\}.$$ Take $\lambda=\rho$, then the polytope cut out is given by $$c_1\geq 0, c_2\geq c_3 \geq 0,1+c_3 \geq c_2,1+c_2 \geq c_1+2c_3,1\geq c_3.$$ which is the corresponding string polytope. Note also that this polytope has integral points $$\{(0, 0, 0), (0, 1, 0), (0, 2, 1), (0, 1, 1), (1, 1, 0), (1, 0, 0), (1,
   2, 1), (2, 1, 0)\}.$$\\
\end{example}

The second way to construct a polytope given a choice of $\lambda\in P^+$ and reduced expression $\textbf{i}$ for $w_0$ is using the theory of Newton-Okounkov bodies [Ok96, Ok98]. Given $\lambda\in P^+$ we can define a line bundle on $G/B$ by $$\mathcal{L}_{\lambda}:=G \times^{B} \mathbb{C}_{-\lambda}=\{(g,z)\}/(g,z)\sim (gb,\lambda(b)z)\qquad \text{  (see [Sp])}$$ and recall that  $H^{0}(G/B,\mathcal{L}_{\lambda})^*=V_{\lambda}$ is the irreducible representation of $G$ with highest weight $\lambda$. Now let $\textbf{i}$ be a  reduced expression for $w_0$ and define $$\nu_{\lambda,\textbf{i}}:H^{0}(G/B,\mathcal{L}_{\lambda})\rightarrow \mathbb{R}^N$$ via the composition\begin{equation}H^{0}(G/B,\mathcal{L}_{\lambda})\xrightarrow{1/\sigma_{\mathrm{lw}}} \mathbb{C}(G/B) \xrightarrow{(\textbf{y}_{\textbf{i}})^*}  \mathbb{C}[x_1^{\pm},...,x_N^{\pm}]\xrightarrow{val} \mathbb{Z}^N \subset \mathbb{R}^N\end{equation} The first map takes a section $\sigma$ to the rational function $\frac{\sigma}{\sigma_{\mathrm{lw}}}$ where $\sigma_{\mathrm{lw}}$ is a lowest weight section (defined up to scalar), $(\textbf{y}_{\textbf{i}})^*$ pulls back this rational function to $(\mathbb{C}^*)^N$ and the map $val$ is the map which selects the exponent of the lexicographically maximal term with respect to the ordering $x_{1}>x_{2}>...>x_{N}$. Let $Q_{\textbf{i}}(\lambda)$ be the image of $\nu_{\lambda,\textbf{i}}$.

The Newton-Okounkov body for this choice of valuation is defined to be $$\mathrm{NO}_{\bf{i}}(\lambda):=\overline{\bigcup_{n=1}^{\infty}\frac{1}{n}Q_{\bf{i}}(n\lambda)}$$ The Newton-Okounkov body $\mathrm{NO}_{\bf{i}}(\lambda)$ is a polytope in $\mathbb{R}^N$ whose integral points are given by $Q_{\bf{i}}(\lambda)$. We then have the following result of Kaveh [Ka].
\begin{theorem}
The Newton-Okounkov body for this choice of valuation $\mathrm{NO}_{\bf{i}}(\lambda)$ coincides with the string polytope $\mathrm{String}_{\bf{i}}(\lambda)$.
\end{theorem}

\begin{example}
Let $G=SL_3$, $\textbf{i}=(212)$ and $\lambda=\rho$. Then $V_{\lambda}$ is the adjoint representation $\mathfrak{g}$ which has a basis given by $T_{1}=diag(1,-1,0),T_{2}=diag(0,1,-1)$ and the root space generators: $$V_{\lambda}=\langle E_{\lambda},E_{\alpha_{1}},E_{\alpha_{2}},T_{1},T_{2},F_{\alpha_{1}},F_{\alpha_{2}},F_{\lambda}\rangle_{\mathbb{C}}.$$ Thus we have the dual basis for $V_{\lambda}^{*}$:$$V_{\lambda}=\langle E_{\lambda}^{*},E_{\alpha_{1}}^{*},E_{\alpha_{2}}^{*},T_{1}^{*},T_{2}^{*},F_{\alpha_{1}}^{*},F_{\alpha_{2}}^{*},F_{\lambda}^{*}\rangle_{\mathbb{C}}.$$ Now $\mathcal{L}_{\lambda}$ is the pullback of $\mathcal{O}(1)$ under the map $G/B \hookrightarrow \mathbb{P}(V_{\lambda})$ which sends $gB \mapsto \langle g \cdot v_{\lambda}^{+} \rangle$, where $v_{\lambda}^{+}$ is a highest weight vector of $V_{\lambda}$. Consider the map $(\mathbb{C}^*)^3\overset{\textbf{y}_{\textbf{i}}}\longrightarrow G/B \hookrightarrow \mathbb{P}(V_{\lambda})$. This map sends \begin{align*}
(x_{1},x_{2},x_{3})&\mapsto \begin{pmatrix}1&0&0\\x_{2}&1&0\\x_{1}x_{2}&x_{1}+x_{3}&1\end{pmatrix}B\\
& \mapsto \left < \begin{pmatrix}1&0&0\\x_{2}&1&0\\x_{1}x_{2}&x_{1}+x_{3}&1\end{pmatrix} \begin{pmatrix}0&0&1\\0&0&0\\0&0&0\end{pmatrix} \begin{pmatrix}1&0&0\\x_{2}&1&0\\x_{1}x_{2}&x_{1}+x_{3}&1\end{pmatrix}^{-1} \right > \\
&=\left < \begin{pmatrix} x_{2}x_{3}&-(x_{1}+x_{3})&1\\x_{2}^{2}x_{3}&-x_{2}x_{3}-x_{1}x_{2}&x_{2}\\x_{1}x_{2}^{2}x_{3}&-x_{1}x_{2}(x_{1}+x_{3})&x_{1}x_{2}\end{pmatrix}\right >.
\end{align*}
From this we can compute the image of the basis of $V_{\lambda}^{*}$ under $\nu_{\lambda,\textbf{i}}$ to be $$\{(0, 0, 0), (0, 1, 0), (0, 2, 1), (0, 1, 1), (1, 1, 0), (1, 0, 0), (1,
   2, 1), (2, 1, 0)\}$$which has convex hull given by the corresponding string polytope.\\
\end{example}

\section{Critical points of $\mathcal{W}$}

Consider $\mathcal{W}_{t^{\lambda}}:Z_{t^{\lambda}}(\mathcal{K}_{>0})\rightarrow \mathcal{K}_{>0}$ and define $Z_{t^{\lambda}}(\mathcal{K}_{>0})^{crit}$ to be the set of critical points of $\mathcal{W}_{t^{\lambda}}$.
The following is the main theorem of this section.
\begin{theorem}\label{tcp}
Let $\lambda\in P^+$. The subset $Z_{t^{\lambda}}(\mathcal{K}_{>0})^{crit}$ consists of a single point, which we denote $p_{\lambda}$. Furthermore, $p_{\lambda}\in Z_{t^{\lambda}}(\mathcal{K}_{>0})^+$ and $\mathrm{wt}(p_{\lambda})=\rm{Id}$.
\end{theorem}
Suppose $\theta:S\rightarrow Z$ is a chart in $\Theta$. Let $$S_{t^{\lambda}}(\mathcal{K}_{>0}):=\{s\in S(\mathcal{K}_{>0}) :\rm{hw}\circ\theta(s)=t^{\lambda}\}$$ and let $S_{t^{\lambda}}(\mathcal{K}_{>0})^{crit}$ be the critical points of $\mathcal{W}_{t^{\lambda}}\circ \theta:S_{t^{\lambda}}(\mathcal{K}_{>0})\rightarrow \mathcal{K}_{>0}$. To prove the theorem it is sufficient to show that $S_{t^{\lambda}}(\mathcal{K}_{>0})^{crit}$ consists of a unique critical point for some chart $\theta:S\rightarrow Z$ in $\Theta$. This is because $\theta$ given a bijection between $S_{t^{\lambda}}(\mathcal{K}_{>0})^{crit}$ and $Z_{t^{\lambda}}(\mathcal{K}_{>0})^{crit}$. With this in mind, we will define, following [R], a toric chart in $\Theta_{Z}$ for which the maps $\mathcal{W}$, $\mathrm{hw}$ and $\mathrm{wt}$ have a particularly nice form.\\

\noindent\textbf{A special toric chart on $Z$}\\
Consider the quiver, first introduced in [G], with $n(n+1)/2$ vertices in lower triangular form with arrows going up and left. Label the vertices as $v_{ij}$ for $1 \leq j \leq i \leq n$ where the labelling is as in matrix entries. Let $\mathcal{V}^{*}=\{v_{11},...,v_{nn}\}$, $\mathcal{V}^{\bullet}=\{v_{ij}:1 \leq j < i \leq n\}$ and $\mathcal{V}=\mathcal{V}^{*} \sqcup \mathcal{V}^{\bullet}$. Let $\mathcal{A}$ be the set of arrows of the quiver. Given an arrow $a\in\mathcal{A} $ denote by $\rm{h}(a)$ and $\rm{t}(a) \in \mathcal{V}$ the head and tail of $a$. When $n=4$ the quiver looks like:

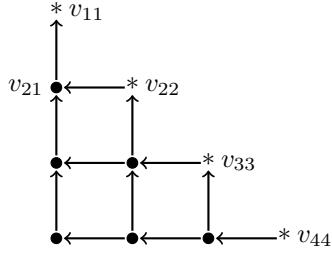
\begin{center}
 \begin{tikzpicture}
\draw [fill] (0.5,0.5) circle [radius=0.07];
\draw [fill] (0.5,1.5) circle [radius=0.07];
\draw [fill] (0.5,2.5) circle [radius=0.07];
\draw [fill] (1.5,0.5) circle [radius=0.07];
\draw [fill] (1.5,1.5) circle [radius=0.07];
\draw [fill] (2.5,0.5) circle [radius=0.07];
\node  at    (0.5,3.5) {*} ;
\node  at    (1.5,2.5) {*}  ;
\node  at    (2.5,1.5) {*}  ;
\node  at    (3.5,0.5) {*}  ;
\node  at    (0.9,3.5) {$v_{11}$} ;
\node  at    (1.9,2.5) {$v_{22}$} ;
\node  at    (2.9,1.5) {$v_{33}$} ;
\node  at    (3.9,0.5) {$v_{44}$} ;
\node  at    (0.1,2.5) {$v_{21}$} ;
\draw[->][thick,black] (0.5,0.6) -- (0.5,1.4);
\draw[->][thick,black] (0.5,1.6) -- (0.5,2.4);
\draw[->][thick,black] (0.5,2.6) -- (0.5,3.4);
\draw[->][thick,black] (1.5,0.6) -- (1.5,1.4);
\draw[->][thick,black] (1.5,1.6) -- (1.5,2.4);
\draw[->][thick,black] (2.5,0.6) -- (2.5,1.4);
\draw[->][thick,black] (3.4,0.5) -- (2.6,0.5);
\draw[->][thick,black] (2.4,0.5) -- (1.6,0.5);
\draw[->][thick,black] (1.4,0.5) -- (0.6,0.5);
\draw[->][thick,black] (2.4,1.5) -- (1.6,1.5);
\draw[->][thick,black] (1.4,1.5) -- (0.6,1.5);
\draw[->][thick,black] (1.4,2.5) -- (0.6,2.5);
\end{tikzpicture}
\end{center}

Consider the torus $(\mathcal{K}^{*})^{\mathcal{V}}$ with coordinates $x_{v}$ for $v \in \mathcal{V}$ and also the torus $$\mathcal{M}:=\{(z_{a})_{a\in \mathcal{A}}\in(\mathcal{K}^{*})^{\mathcal{A}}:z_{a_{1}}z_{a_{2}}=z_{a_{3}}z_{a_{4}} \text{ when } a_{1},a_{2},a_{3},a_{4} \text{ form a square as in Fig.1 }\}
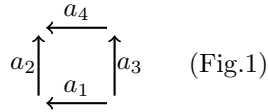
$$
\begin{center}
 \begin{tikzpicture}
\node  at    (3,1.0)  {(Fig.1)};
\node  at    (0.3,1) {$a_{2}$} ;
\node  at    (1.7,1) {$a_{3}$}  ;
\node  at    (1,1.7) {$a_{4}$}  ;
\node  at    (1,0.7) {$a_{1}$}  ;
\draw[->][thick,black] (0.5,0.6) -- (0.5,1.4);
\draw[->][thick,black] (1.5,0.6) -- (1.5,1.4);
\draw[->][thick,black] (1.4,0.5) -- (0.6,0.5);
\draw[->][thick,black] (1.4,1.5) -- (0.6,1.5);
\end{tikzpicture}
\end{center}
Note we have a projection $(\mathcal{K}^{*})^{\mathcal{V}}\rightarrow \mathcal{M}$ defined by $z_{a}=\frac{x_{h(a)} }{ x_{t(a)}}$. Next define three maps $$\kappa:\mathcal{M}\rightarrow T^{\vee}(\mathcal{K}) \qquad \mathcal{F}:\mathcal{M}\rightarrow\mathcal{K} \qquad \gamma:\mathcal{M}\rightarrow T^{\vee}(\mathcal{K})$$
\text{\tiny $\bullet$} The map $\kappa$ is defined to be the map induced by the map $$\tilde{\kappa}:(\mathcal{K}^{*})^{\mathcal{V}}\rightarrow T_{GL_{n}^{\vee}}(\mathcal{K}) \quad\text{   sending   } \quad(x_{v})_{v \in \mathcal{V}}\mapsto (x_{v_{ii}})_{i=1,..,n}$$ Note this induced map is well-defined since $T^{\vee}\subset PGL_{n}$.\\

\noindent \text{\tiny $\bullet$}  The map $\mathcal{F}$ is defined by $(z_{a})_{a \in \mathcal{A}}\mapsto \sum\limits_{a \in \mathcal{A}} z_{a}$\\

\noindent\text{\tiny $\bullet$}  The map $\gamma$ is defined as follows: For $i=1,..,n$ let $\mathcal{D}_{i}=\{v_{i,1},v_{i+1,2},...,v_{n,n+1-i}\}$ be the $i$th diagonal and let $\zeta_{i}=\prod\limits_{v \in \mathcal{D}_{i}}x_{v}$. We set $\zeta_{n+1}=1$. Then we define a map $$\tilde{\gamma}:(\mathcal{K}^{*})^{\mathcal{V}}\rightarrow T_{GL_{n}^{\vee}}(\mathcal{K})\quad \text{   sending   } \quad(x_{v})_{v \in \mathcal{V}}\mapsto (t_{i})_{i=1,..,n}\quad \text{   where   }\quad t_{i}=\frac{\zeta_{i}}{\zeta_{i+1}}$$ and note that this map descends to give a well defined map $\mathcal{M}\rightarrow T^{\vee}(\mathcal{K})$.\\

Now we will define a toric chart $\theta_{\mathcal{M}}:\mathcal{M}\rightarrow Z(\mathcal{K})$. First some other definitions. For a reduced expression $\textbf{i}=(i_{1},...,i_{N})$ of $w_0$, let $$\textbf{x}_{\textbf{i}}:(\mathcal{K}^*)^N \rightarrow (U^{\vee}\cap B_{-}^{\vee}\bar{w}_0B_{-}^{\vee})(\mathcal{K}) \quad \text{   be the map   } \quad(y_{1},...,y_{N}) \mapsto \textbf{x}_{i_{1}}(y_{1})...\textbf{x}_{i_{N}}(y_{N})$$

Consider the particular choice of reduced expression given by $$\textbf{i}_{0}=(1,2,..,n-1,1,2,..,n-2,.....,1,2,1)$$
Also let $z_{ij}$ be the coordinate $z_{a}$ where $a$ is the vertical arrow with $h(a)=v_{ij}$.

Now define a toric chart $\theta_{\mathcal{M}}:\mathcal{M}\rightarrow Z(\mathcal{K})$ given by  $$z=(z_{a})_{a \in \mathcal{A}} \mapsto  \Phi(\kappa(z),\textbf{x}_{\textbf{i}_0}(z_{n-1 \, 1},z_{n-2 \,1},..,z_{1\, 1},z_{n-1\, 2},z_{n-2\, 2},..,z_{2 \, 2},....,z_{n-1\, n-2},z_{n-2\, n-2},z_{n-1\, n-1}))$$

\begin{remark}
This sequence of arrows $z_{n-1 \, 1},z_{n-2 \,1},..$ can be described by the sequence of vertical arrows starting at the bottom left corner and moving vertically upward until you reach the top and then moving the bottom of the next column and repeating.
\end{remark}

The toric chart $\theta_{\mathcal{M}}$ is positively compatible with $\Theta$, hence lies in the positive structure $\Theta$. This again follows from [BZ]. The reason this toric chart is particularly nice because of the following Lemma, which follows from [R,Thrm 9.2] and [R,Thrm 9.7].
\begin{lemma}With the above notation, we have: $$\mathcal{W} \circ \theta_{\mathcal{M}}=\mathcal{F}\quad\text{and}\quad  \mathrm{hw} \circ \theta_{\mathcal{M}}=\kappa  \quad\text{and}\quad  \mathrm{wt} \circ \theta_{\mathcal{M}}=\gamma$$
\end{lemma}
\noindent\textit{Proof of Lemma 5.3}:
Let $$\mathcal{W}':Z(\mathcal{K})\rightarrow \mathcal{K}\quad \text{sending}\quad u_1q\bar{w}_{0}u_2\mapsto \chi(u_1)+\chi(u_2)$$ and $$\mathrm{hw}':Z(\mathcal{K})\rightarrow T^{\vee}(\mathcal{K})\quad \text{sending}\quad u_1q\bar{w}_{0}u_2\mapsto q$$ and $$\mathrm{wt}':Z(\mathcal{K})\rightarrow T^{\vee}(\mathcal{K})\quad \text{sending}\quad b\mapsto pr(b)$$ Also let $$\textbf{i}_{0}'=(n-1,n-2,..,1,n-1,n-2,..,2,.....,n-1,n-2,n-1)$$ and define $\theta_{\mathcal{M}}':\mathcal{M}\rightarrow Z(\mathcal{K})$ to be given by  $$z=(z_{a})_{a \in \mathcal{A}} \mapsto  \Phi(\kappa(z),\textbf{x}_{\textbf{i}_{0}'}(z_{n-1 \, 1},z_{n-2 \,1},..,z_{1\, 1},z_{n-1\, 2},z_{n-2\, 2},..,z_{2 \, 2},....,z_{n-1\, n-2},z_{n-2\, n-2},z_{n-1\, n-1}))$$ It shown in [R,Thrm 9.2] that $\mathcal{W}' \circ \theta_{\mathcal{M}'}=\mathcal{F}$.
Now define a map $$\iota:G^{\vee}\rightarrow G^{\vee}\quad \text{sending}\quad g\mapsto (\bar{w}_0 g^{-1} \bar{w}_0^{-1})^T$$
where $()^T$ is the transpose. The map $\iota$ acts on $U^{\vee}$ and also on $Z$ and observe that for $u\in U^{\vee}$ we have that $\chi(\iota(u))=\chi(u)$, where $\chi$ is the map from Section 2. Multiplying out we can see that $$\iota(u_1q\bar{w}_{0}u_2)=\iota(u_1)\bar{w}_{0}q^{-1}\iota(u_2)$$ hence $\mathcal{W}'=\mathcal{W}\circ \iota$.
Also, we have that $$\theta_{\mathcal{M}}=\iota\circ\theta_{\mathcal{M}}'$$ since $\iota(\textbf{x}_{\textbf{i}_{0}'}(z_{n-1 \, 1},z_{n-2 \,1},...,z_{n-2\, n-2},z_{n-1\, n-1}))=\textbf{x}_{\textbf{i}_{0}}(z_{n-1 \, 1},z_{n-2 \,1},...,z_{n-2\, n-2},z_{n-1\, n-1})$. So, combining these facts, we get $\mathcal{W}\circ\theta_{\mathcal{M}}= \mathcal{W}\circ\iota\circ\theta_{\mathcal{M}}'=\mathcal{W}'\circ\theta_{\mathcal{M}}'=\mathcal{F}$.

The identity $\mathrm{wt} \circ \theta_{\mathcal{M}}=\gamma$ is proved similarly, using the fact $\mathrm{wt}' \circ \theta_{\mathcal{M}}'=\gamma$ which follows from [R,Thrm 9.7]. The identity $\mathrm{hw} \circ \theta_{\mathcal{M}}=\kappa$ follows immediately from the definitions. $\square$

\begin{example}
Let $n=3$, then $\textbf{i}_{0}=(121)$ and the quiver looks like:

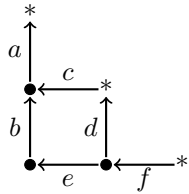
\begin{center}
 \begin{tikzpicture}
\draw [fill] (0.5,0.5) circle [radius=0.07];
\draw [fill] (0.5,1.5) circle [radius=0.07];
\draw [fill] (1.5,0.5) circle [radius=0.07];
\node  at    (0.5,2.5) {*} ;
\node  at    (1.5,1.5) {*}  ;
\node  at    (2.5,0.5) {*}  ;
\node[fill=white] at    (0.3,1) {$b$}  ;
\node[fill=white] at    (0.3,2) {$a$}  ;
\node[fill=white] at    (1.3,1) {$d$}  ;
\node[fill=white] at    (1,0.3) {$e$}  ;
\node[fill=white] at    (1,1.7) {$c$}  ;
\node[fill=white] at    (2,0.3) {$f$}  ;
\draw[->][thick,black] (0.5,0.6) -- (0.5,1.4);
\draw[->][thick,black] (0.5,1.6) -- (0.5,2.4);
\draw[->][thick,black] (1.5,0.6) -- (1.5,1.4);
\draw[->][thick,black] (2.4,0.5) -- (1.6,0.5);
\draw[->][thick,black] (1.4,0.5) -- (0.6,0.5);
\draw[->][thick,black] (1.4,1.5) -- (0.6,1.5);
\end{tikzpicture}
\end{center}
Then $$\theta_{\mathcal{M}}:(a,b,c,d,e,f)\mapsto \begin{pmatrix}\frac{1}{ef}&&\\ \frac{1}{cdf}&\frac{1}{bf}&\\ \frac{1}{acdf}&\frac{b+d}{abdf}&\frac{1}{ad} \end{pmatrix}$$
and we can compute $\mathrm{hw},\mathcal{W},\mathrm{wt}$ in these coordinates (using Example 2.2):
$$(a,b,c,d,e,f)\mapsto \begin{pmatrix}acdf&&\\ &df&\\ &&1 \end{pmatrix} $$
$$(a,b,c,d,e,f)\mapsto a+b+c+d+e+f $$
$$(a,b,c,d,e,f)\mapsto \begin{pmatrix}ad&&\\ &bf&\\ &&fe\end{pmatrix} $$
\end{example}

So to prove Theorem 5.1 we need to show that $\mathcal{F}:\mathcal{M}_{t^{\lambda}}(\mathcal{K}_{>0})\rightarrow\mathcal{K}_{>0}$ has a unique critical point. In fact, it will be more convenient to work with the vertex coordinates. Consider the maps $\tilde{\kappa}:(\mathcal{K}^{*})^{\mathcal{V}}\rightarrow T_{GL_{n}^{\vee}}(\mathcal{K})$ and $\tilde{\gamma}:(\mathcal{K}^{*})^{\mathcal{V}}\rightarrow T_{GL_{n}^{\vee}}(\mathcal{K})$, and think of $\mathcal{F}$ as a function from $(\mathcal{K}^*)^{\mathcal{V}}$ to $\mathcal{K}$. Let $\tilde{\lambda}=(\lambda_{1},\lambda_{2},..,\lambda_{n})\in X_{*}(T_{GL_{n}^{\vee}})$ be a lift of $\lambda$. We will show that the fiber of $(\mathcal{K}^*)^{\mathcal{V}}$ lying over $\tilde{\lambda}$ contains a unique totally positive critical point. This will prove what we want.

The condition for $(x_v)_{v\in \mathcal{V}}$ to lie in the fiber over $t^{\tilde{\lambda}}$ is given by $x_{v_{ii}}=t^{\lambda_{i}}$ for $i=1,...,n$. Now $x_v$ for $v\in\mathcal{V}^{\bullet}$ give a system of coordinates on the fiber over $\tilde{\lambda}$, so we can use them to compute the critical points. We have $$x_v\frac{\partial \mathcal{F}}{\partial x_v}=\sum_{a\in\mathcal{A}:\mathrm{h}(a)=v}\frac{x_{\mathrm{h}(a)}}{x_{\mathrm{t}(a)}}\quad-\sum_{a\in\mathcal{A}:\mathrm{t}(a)=v}\frac{x_{\mathrm{h}(a)}}{x_{\mathrm{t}(a)}}$$
which give the \textit{critical point conditions}:
\begin{equation}\label{e:cpconds}
\sum_{a:\mathrm{h}(a)=v}\frac{x_{\mathrm{h}(a)}}{x_{\mathrm{t}(a)}}\quad=\sum_{a:\mathrm{t}(a)=v}\frac{x_{\mathrm{h}(a)}}{x_{\mathrm{t}(a)}}\quad \text{ for } v\in\mathcal{V}^{\bullet}
\end{equation}
So we need to prove that the set $$\mathcal{C}_{\lambda}:=\{ (x_v)\in (\mathcal{K}_{>0})^{\mathcal{V}}: x_{v_{ii}}=t^{\lambda_{i}} \text{ for } v=v_{ii}\in\mathcal{V}^* \quad \text{and} \sum_{a:\mathrm{h}(a)=v}\frac{x_{\mathrm{h}(a)}}{x_{\mathrm{t}(a)}}\quad=\sum_{a:\mathrm{t}(a)=v}\frac{x_{\mathrm{h}(a)}}{x_{\mathrm{t}(a)}} \text{ for } v\in\mathcal{V}^{\bullet} \}$$
consists of a single point. This will involve two parts, existence and uniqueness.

\subsection*{Existence}
First we will inductively define $\delta_v\in\mathbb{R}$ for $v\in \mathcal{V}$. Let $\delta_v=\lambda_{i}$ for $v=v_{ii}\in\mathcal{V}^{*}$. Let $\mathcal{V}_{0}^{\bullet}=\mathcal{V}^{*}$ and $\mathcal{A}_{0}^{\bullet}=\emptyset$. Let $\mathcal{V}_{0}=\mathcal{V}_{0}^{\bullet}$ and $\mathcal{A}_{0}=\mathcal{A}_{0}^{\bullet}$.

Given a path $\pi$ which follows the directed arrows of the quiver, let $\mathrm{ver}(\pi)$ be the set of vertices contained in the path and let $\mathrm{s}(\pi)$ be the vertex where $\pi$ starts and $\mathrm{e}(\pi)$ be the vertex where $\pi$ ends. Let $\mathrm{arr}(\pi)$ be the set of arrows making up $\pi$ and $\mathrm{len}(\pi)$ be the length of the path, i.e number of arrows in it.

For $l \geq 1$ let $$\Gamma_{l}=\{\text{paths } \pi:\mathrm{s}(\pi)\in \mathcal{V}_{l-1} \text{ and } \mathrm{e}(\pi)\in \mathcal{V}_{l-1} \text{ and }\mathrm{arr}(\pi) \cap  \mathcal{A}_{l-1}=\emptyset \}$$
Let $$\gamma_{l}:\Gamma_{l}\rightarrow \mathbb{Q}\text{ be the map }\pi \mapsto \frac{\delta_{e(\pi)}-\delta_{s(\pi)}}{\mathrm{len}(\pi)}$$
and $$\kappa_{l}=\min_{\pi\in\Gamma_{l}}\gamma_{l}$$
Let $$\mathcal{A}_{l}^{\bullet}=\bigcup_{\pi \in \Gamma_{l}:\gamma_{l}(\pi)=\kappa_{l}}\mathrm{arr}(\pi)$$
and $$\mathcal{V}_{l}^{\bullet}=\bigcup_{\pi \in \Gamma_{l}:\gamma_{l}(\pi)=\kappa_{l}}\mathrm{ver}(\pi)\backslash\mathcal{V}_{l-1}$$

For $a\in\mathcal{A}_{l}^{\bullet}$ let $\sigma_a=\kappa_l$ and for $v \in \mathcal{V}_{l}^{\bullet}$ define $\delta_v$ via $\sigma_a=\delta_{\mathrm{h}(a)}-\delta_{\mathrm{t}(a)}$. Also let $\mathcal{A}_{l}=\mathcal{A}_{l-1}\sqcup\mathcal{A}_{l}^{\bullet}$ and $\mathcal{V}_{l}=\mathcal{V}_{l-1}\sqcup\mathcal{V}_{l}^{\bullet}$.

For some $l$ we will have $\mathcal{A}_l=\mathcal{A}$, at which point we will have defined $\delta_v$ for all $v\in \mathcal{V}$. We also set $\sigma_a=\delta_{\mathrm{h}(a)}-\delta_{\mathrm{t}(a)}$ for $a \in \mathcal{A}$. Observe that $\kappa_{l+1}>\kappa_l$, hence if $v\in\mathcal{V}_{l}^{\bullet}$ then we have $\min_{a\in\mathcal{A}:\mathrm{h}(a)=v}\sigma_a\quad=\kappa_l=\min_{a\in\mathcal{A}:\mathrm{t}(a)=v}\sigma_a$. For $v\in\mathcal{V}^{\bullet}$ define $$\pi(v):=\min_{a\in\mathcal{A}:\mathrm{h}(a)=v}\sigma_a=\min_{a\in\mathcal{A}:\mathrm{t}(a)=v}\sigma_a$$
We say these $\sigma_a$ satisfy the ``\textit{tropical critical point conditions}": \begin{equation}\label{e:tcpconds}\min_{a\in\mathcal{A}:\mathrm{h}(a)=v}\sigma_a\quad=\min_{a\in\mathcal{A}:\mathrm{t}(a)=v}\sigma_a\quad \text{for } v\in \mathcal{V}^{\bullet}\end{equation}

Let $M$ be the lcm of the demoninators of the rational numbers $\sigma_a$ (in lowest terms). Now, for $v\in \mathcal{V}$, let $$x_v=d_v t^{\delta_v}\sum\limits_{k\geq 0}x_{v,k}t^{k/M}$$ where the $d_v$ and $x_{v,k}$ are for the moment just variables. For all $v\in \mathcal{V}$ set $x_{v,0}=1$ and for $v \in \mathcal{V}^*$ set $x_{v,k}=0$ for $k\geq1$.

Let $z_a=\frac{x_{\mathrm{h}(a)}}{x_{\mathrm{t}(a)}}$, then $z_a=c_a t^{\sigma_a}\sum\limits_{k \geq 0}z_{a,k}t^{k/M}$ where
\begin{align*}
&z_{a,0}=1\quad \text{   and }\\&z_{a,k}=x_{\mathrm{h}(a),k}-x_{\mathrm{t}(a),k}-\sum\limits_{j=1}^{k-1}x_{\mathrm{t}(a),j}z_{a,k-j}\quad\text{   for  }k\geq1
\end{align*}

Now consider $$\mathcal{F}_{\mathcal{A}_{l}^{\bullet}}:\mathbb{R}_{>0}^{\mathcal{V}_{l}^{\bullet}}\rightarrow\mathbb{R}_{>0} \text{ sending }(d_v)\mapsto\sum_{a\in\mathcal{A}_{l}^{\bullet}}\frac{d_{\mathrm{h}(a)}}{d_{\mathrm{t}(a)}}$$ where we think of $d_v$ for $v\in\mathcal{V}_{l-1}$ as already determined constants in $\mathbb{R}_{>0}$. Also let $c_{a}=\frac{d_{\mathrm{h}(a)}}{d_{\mathrm{t}(a)}}$.

For $v\in\mathcal{V}^{\bullet}_{l}$ we have $$d_v\frac{\partial \mathcal{F}_{\mathcal{A}^{\bullet}_{l}}}{\partial d_v}=\sum_{a\in\mathcal{A}^{\bullet}_{l}:\mathrm{h}(a)=v}\frac{d_{\mathrm{h}(a)}}{d_{\mathrm{t}(a)}}\quad-\sum_{a\in\mathcal{A}^{\bullet}_{l}:\mathrm{t}(a)=v}\frac{d_{\mathrm{h}(a)}}{d_{\mathrm{t}(a)}}$$
These give the critical point conditions for $ \mathcal{F}_{\mathcal{A}^{\bullet}_{l}}$.

\begin{lemma}
$\mathcal{F}_{\mathcal{A}_{l}^{\bullet}}$ has a unique critical point in $\mathbb{R}_{>0}^{\mathcal{V}_{l}^{\bullet}}$.
\end{lemma}
\noindent\textit{Proof of Lemma 5.5}:
This statement is essentially proven in [R06,Theorem 10.2], and here we follow the argument from that proof. To show existence, define a sequence of compact subsets $S_1\subset S_2 \subset ...\subset\mathbb{R}_{>0}^{\mathcal{V}_{l}^{\bullet}} $ by $$S_{m}:=\{(d_v)\in\mathbb{R}_{>0}^{\mathcal{V}_{l}^{\bullet}}: m^{-\vert \mathcal{A}_{l}^{\bullet} \vert}\min\limits_{ w\in\mathcal{V}_{l-1}}d_w\,\leq\, d_v \,\leq \, m^{\vert \mathcal{A}_{l}^{\bullet} \vert} \max\limits_{ w\in\mathcal{V}_{l-1}}d_w \text{   for all $v \in \mathcal{V}_{l}^{\bullet}$ }\}$$ If $(d_v)\notin S_m$ then there is some vertex $v$ for which one of the inequalities fails to hold. Let $\pi\in \Gamma_{l}$ be some path with $v \in \mathrm{ver}(\pi)$. If the first inequality fails to hold then the product of $c_a$ for the arrows in $\pi$ from $v$ to $\mathrm{e}(\pi)$ equals $\frac{d_{\mathrm{e}(\pi)}}{d_v}>\frac{d_{\mathrm{e}(\pi)}m^{\vert \mathcal{A}_{l}^{\bullet} \vert}}{\min\limits_{ w\in\mathcal{V}_{l-1}}d_w} \geq m^{\vert \mathcal{A}_{l}^{\bullet} \vert}$, so there is some arrow in $\pi$ such that $c_a>m$. If the other inequality fails to hold then the product of $c_a$ for the arrows in $\pi$ from $\mathrm{s}(\pi)$ to $v$ equals $\frac{d_v}{d_{\mathrm{s}(\pi)}}>\frac{ m^{\vert \mathcal{A}_{l}^{\bullet} \vert} \max\limits_{ w\in\mathcal{V}_{l-1}}d_w}{d_{\mathrm{s}(\pi)}}\geq m^{\vert \mathcal{A}_{l}^{\bullet} \vert} $, so again there is some arrow in $\pi$ such that $c_a>m$. Hence $(d_v)\notin S_m$ implies $\mathcal{F}_{\mathcal{A}_{l}^{\bullet}}>m$. Now for every $m$ with $S_m$ non-empty, $\mathcal{F}_{\mathcal{A}_{l}^{\bullet}}$ attains a minimum $p_m$, and the sequence of minima $p_m\geq p_{m+1}\geq...$ stabilises to a global minimum, since eventually $p_m<m$.

For uniqueness it is sufficient to show that the Hessian of $\mathcal{F}_{\mathcal{A}^{\bullet}_{l}}$ is everywhere positive definite, which follows by direct calculation:
\begin{align*}
\left ( \sum\limits_{v\in\mathcal{V}^{\bullet}_{l}}m_v\frac{\partial}{\partial d_v}\right )^2 \mathcal{F}_{\mathcal{A}^{\bullet}_{l}} &=  \left ( \sum\limits_{v\in\mathcal{V}^{\bullet}_{l}}m_v\frac{\partial}{\partial d_v}\right )\sum\limits_{v'\in\mathcal{V}^{\bullet}_{l}}m_{v'}\left (  \sum_{a\in\mathcal{A}^{\bullet}_{l}:\mathrm{h}(a)=v'}\frac{d_{\mathrm{h}(a)}}{d_{\mathrm{t}(a)}}-
\sum_{a\in\mathcal{A}^{\bullet}_{l}:\mathrm{t}(a)=v'}\frac{d_{\mathrm{h}(a)}}{d_{\mathrm{t}(a)}}
  \right)\\
&= \sum\limits_{v\in\mathcal{V}^{\bullet}_{l}}m_v^2\left (  \sum_{a\in\mathcal{A}^{\bullet}_{l}:\mathrm{h}(a)=v}\frac{d_{\mathrm{h}(a)}}{d_{\mathrm{t}(a)}}+
\sum_{a\in\mathcal{A}^{\bullet}_{l}:\mathrm{t}(a)=v}\frac{d_{\mathrm{h}(a)}}{d_{\mathrm{t}(a)}}
  \right)-2\sum\limits_{a\in\mathcal{A}_{l}^{\bullet}}m_{\mathrm{h}(a)}m_{\mathrm{t}(a)}\frac{d_{\mathrm{h}(a)}}{d_{\mathrm{t}(a)}}\\
&=\sum\limits_{a\in\mathcal{A}_{l}^{\bullet}}(m_{\mathrm{h}(a)}-m_{\mathrm{t}(a)})^2\frac{d_{\mathrm{h}(a)}}{d_{\mathrm{t}(a)}}\quad
\quad\quad\quad\quad\quad\quad\quad\quad\quad\quad\quad\quad\quad\quad\quad\quad\quad\square
\end{align*}\\

Now for $v\in \mathcal{V}^{\bullet}$ define $$\mathrm{crit}(v):= \sum\limits_{\substack{a\in\mathcal{A}\\ \mathrm{h}(a)=v}}  \frac{x_{\mathrm{h}(a)}}{x_{\mathrm{t}(a)}} - \sum\limits_{\substack{a\in\mathcal{A}\\ \mathrm{t}(a)=v}}  \frac{x_{\mathrm{h}(a)}}{x_{\mathrm{t}(a)}}$$ and $$\mathrm{coeff}(v,k):=\text{ the coefficient of }t^{\pi(v)+k/M}\text{ in }\mathrm{crit}(v)$$\\

Next we inductively assign $d_v\in\mathbb{R}_{>0}$ for $v\in \mathcal{V}$. Let $d_v=1$ for $v\in\mathcal{V}^{*}$. Then we have the following claim: Suppose we assigned values in $\mathbb{R}_{>0}$ to $d_v$ for $v\in\mathcal{V}_{l-1}$ such that $\mathrm{coeff}(v,0)=0$ for all $v\in\mathcal{V}_{l-1}$, then there is a unique choice of $(d_v)_{v\in\mathcal{V}_{l}^{\bullet}}\in (\mathbb{R}_{>0})^{\mathcal{V}_{l}^{\bullet}}$ such that $\mathrm{coeff}(v,0)=0$ for all $v\in\mathcal{V}_{l}$.

This follows from Lemma 5.5, since for $v\in\mathcal{V}^{\bullet}_{l}$ we have
\begin{align*}
\mathrm{coeff}(v,0)&=\sum\limits_{\substack{a:\mathrm{h}(a)=v \\ \sigma_a=\pi(v)}}\frac{d_{\mathrm{h}(a)}}{d_{\mathrm{t}(a)}}-\sum\limits_{\substack{a:\mathrm{t}(a)=v \\ \sigma_a=\pi(v)}}\frac{d_{\mathrm{h}(a)}}{d_{\mathrm{t}(a)}}\\ &=\sum\limits_{a\in\mathcal{A}^{\bullet}_{l}:\mathrm{h}(a)=v}\frac{d_{\mathrm{h}(a)}}{d_{\mathrm{t}(a)}}-\sum\limits_{a\in\mathcal{A}^{\bullet}_{l}:\mathrm{t}(a)=v}\frac{d_{\mathrm{h}(a)}}{d_{\mathrm{t}(a)}}
\end{align*} which gives the critical point equations for $\mathcal{F}_{\mathcal{A}^{\bullet}_{l}}$. So the unique such assignment is given by the unique critical point of  $\mathcal{F}_{\mathcal{A}_{l}}$.

Thus we have assigned values to $\delta_v$ and $d_v$ for all $v\in \mathcal{V}$ such that $\mathrm{coeff}(v,0)=0$ for all $v\in\mathcal{V}^{\bullet}$. We also set $\sigma_a=\delta_{\mathrm{h}(a)}-\delta_{\mathrm{t}(a)}$ and $c_{a}=\frac{d_{\mathrm{h}(a)}}{d_{\mathrm{t}(a)}}$ for $a \in \mathcal{A}$.

Now we have the following claim: Let $k,l\geq 1$ . Suppose we have assigned values to $x_{v,j}$ for $j<k$ and $x_{v,k}$ for $v\in\mathcal{V}^{\bullet}_{m}$ with $m<l$ such that $\mathrm{coeff}(v,j)=0$ for all $v\in \mathcal{V}^{\bullet},j<k$ and $\mathrm{coeff}(v,k)=0$ for all $v\in\mathcal{V}^{\bullet}_{m}$ with $m<l$, then there is a unique way to assign values to $x_{v,k}$ for $v\in\mathcal{V}^{\bullet}_{l}$ such that we also have $\mathrm{coeff}(v,k)=0$ for all $v\in\mathcal{V}^{\bullet}_{l}$. This will complete the proof of existence as it will show that we can inductively define the $x_{v,k}$ so that the critical point conditions hold.

Define a map $$\Psi_{l,k}:\mathbb{R}^{\mathcal{V}_{l}^{\bullet}}\rightarrow\mathbb{R}^{\mathcal{V}_{l}^{\bullet}}\,\,\text{   sending   } \,\, (x_{v,k})_{v\in\mathcal{V}_{l}^{\bullet}}\mapsto (\mathrm{coeff}(v,k))_{v\in\mathcal{V}_{l}^{\bullet}}$$ This is linear in the variables $(x_{v,k})_{v\in\mathcal{V}_{l}^{\bullet}}$. Consider the translated linear map $\Psi_{l,k}^{0}:=\Psi_{l,k}-\Psi_{l,k}(0)$. This map is given by $$\Psi_{l,k}^{0}:\mathbb{R}^{\mathcal{V}_{l}^{\bullet}}\rightarrow\mathbb{R}^{\mathcal{V}_{l}^{\bullet}}\,\,\text{  sending   }\,\, (y_{v})_{v\in\mathcal{V}_{l}^{\bullet}}\mapsto\sum_{\substack{a:\mathrm{h}(a)=v \\ \sigma_a=\pi(v)}}c_a (y_{\mathrm{h}(a)}-y_{\mathrm{t}(a)})\quad-\sum_{\substack{a:\mathrm{t}(a)=v \\ \sigma_a=\pi(v)}}c_a (y_{\mathrm{h}(a)}-y_{\mathrm{t}(a)})$$ where we set $y_{v}=0$ for $v\in\mathcal{V}_{l-1}$. Now we claim this map is positive definite, which can be seen as follows:

\begin{align*}
&\sum\limits_{v\in\mathcal{V}_{l}^{\bullet} }y_v \left (\sum\limits_{\substack{a:\mathrm{h}(a)=v \\ \sigma_a=\pi(v)}}c_a (y_{\mathrm{h}(a)}-y_{\mathrm{t}(a)})\quad-\sum\limits_{\substack{a:\mathrm{t}(a)=v \\ \sigma_a=\pi(v)}}c_a (y_{\mathrm{h}(a)}-y_{\mathrm{t}(a)})\right )
\\&= \sum\limits_{v\in\mathcal{V}_{l}^{\bullet}} \sum\limits_{\substack{a\in\mathcal{A}_{l}^{\bullet}\\ \mathrm{h}(a)=v}}  y_v c_a(y_{\mathrm{h}(a)}-y_{\mathrm{t}(a)})- \sum\limits_{v\in\mathcal{V}_{l}^{\bullet}} \sum\limits_{\substack{a\in\mathcal{A}_{l}^{\bullet}\\ \mathrm{t}(a)=v}}     y_v c_a(y_{\mathrm{h}(a)}-y_{\mathrm{t}(a)})
\\&=\sum\limits_{a\in\mathcal{A}_{l}^{\bullet}} y_{\mathrm{h}(a)}c_a(y_{\mathrm{h}(a)}-y_{\mathrm{t}(a)}) -\sum\limits_{a\in\mathcal{A}_{l}^{\bullet}} y_{\mathrm{t}(a)}c_a(y_{\mathrm{h}(a)}-y_{\mathrm{t}(a)})
\\&=\sum\limits_{a\in\mathcal{A}_{l}^{\bullet}} c_a(y_{\mathrm{h}(a)}-y_{\mathrm{t}(a)})^{2}
\\&\geq0
\end{align*}
with equality iff $y_{\mathrm{h}(a)}=y_{\mathrm{t}(a)}$ for all $a\in\mathcal{A}_{l}^{\bullet}$, which implies all the $y_{v}=0$. Hence $\Psi_{l,k}^{0}$ is invertible, and therefore so is $\Psi_{l,k}$. Thus there is a unique choice of the $x_{v,k}$ for $v\in\mathcal{V}^{\bullet}_{l}$ such that $\mathrm{coeff}(v,l)=0$ for all $v\in\mathcal{V}^{\bullet}_{l}$.

This will be a critical point since after assigning all the coefficients $x_{v,k}$ all the $\mathrm{coeff}(v,j)=0$ for $j\leq k$, and every non-zero term of $\mathrm{crit}(v)$ is of the form $\mathrm{coeff}(v,k)t^{\pi(v)+k/M}$ for some $k\geq0$ since every $\sigma_a=m/M$ for some integer $m$. This completes the proof of existence.

\subsection*{Uniqueness}

For uniqueness, suppose $(x_v)_{\mathcal{V}}\in\mathcal{C}_{\lambda}$. Then $$\sum_{a\in\mathcal{A}:\mathrm{h}(a)=v}\frac{x_{\mathrm{h}(a)}}{x_{\mathrm{t}(a)}}\quad=\sum_{a\in\mathcal{A}:\mathrm{t}(a)=v}\frac{x_{\mathrm{h}(a)}}{x_{\mathrm{t}(a)}}$$ for $v\in \mathcal{V}^{\bullet}$. Let $\delta_v=\mathrm{val}(x_v)$ and $\sigma_a=\delta_{\mathrm{h}(a)}-\delta_{\mathrm{t}(a)}$, so $\delta_v=\lambda_i$ for $v=v_{ii}\in\mathcal{V}^*$. Then applying $\mathrm{val}$ to both sides of this equation we get that the $\sigma_a$ must satisfy the ``tropical critical point conditions'': $$\min_{a\in\mathcal{A}:\mathrm{h}(a)=v}\sigma_a\quad=\min_{a\in\mathcal{A}:\mathrm{t}(a)=v}\sigma_a$$ for $v\in \mathcal{V}^{\bullet}$.

Now we have
\begin{prop}
Given $\lambda\in P^+$, let $\delta_v=\lambda_i$ for $v=v_{ii}\in\mathcal{V}^*$. Then there is a unique point $(\delta_v)\in\mathbb{R}^{\mathcal{V}}$ such that $$\min_{a:\mathrm{h}(a)=v}\sigma_a\quad=\min_{a:\mathrm{t}(a)=v}\sigma_a$$ for $v\in \mathcal{V}^{\bullet}$.
\end{prop}

Given this Proposition, uniqueness follows from the construction of the critical point in the previous section. To see this, suppose $(x_v)\in\mathcal{C}_{\lambda}$, and let $(\delta_v)\in\mathbb{R}^{\mathcal{V}}$ be the unique point satisfying the tropical critical conditions. Then we can write $x_v=d_v t^{\delta_v}(1+\sum\limits_{k\geq 1}x_{v,k}t^{k/M'})$ for some $M'\in\mathbb{Z}_{>0}$, $d_v\in\mathbb{R}_{>0}$ and $x_{v,k}\in\mathbb{R}$. From the previous section on existence, the critical point conditions uniquely specify the values for $d_v$ and $x_{v,k}$, hence there exists a unique critical point.\\

\noindent\textit{Proof of Proposition 5.6}:
The existence of such a point follows from the existence of the critical point, so we just need to show uniqueness. We will prove this by induction on $n$. The case $n=2$ is clear as we must have both arrow coordinates equal to $(\lambda_{1}-\lambda_{2})/2$. Suppose for induction that there is a unique such point for any diagram of size $n-1$, and we will prove it for a diagram of size $n$. Suppose we have two points $(T_{v})_{v \in \mathcal{V}}$ and $(S_{v})_{v \in \mathcal{V}}$ satisfying the tropical critical point conditions with $T_{v_{ii}}=S_{v_{ii}}=\lambda_{i}$ for $i=1,..,n$. Let $\sigma_{a}=T_{h(a)}-T_{t(a)}$ and $\rho_{a}=S_{h(a)}-S_{t(a)}$ and let $\tau_{a}=\sigma_{a}-\rho_{a}$. We want to show $\tau_{a}=0$ for all $a\in \mathcal{A}$.

First we make some definitions. Given two touching arrows $a,b$ define $$\mathrm{sgn}(a,b)=\begin{cases} 0 &\mbox{if } \tau_{a}\tau_{b}=0\\ 1 &\mbox{if } \tau_{a}\tau_{b}>0\\-1&\mbox{if } \tau_{a}\tau_{b}<0 \end{cases}$$ Suppose we have a vertex $v\in\mathcal{V}$
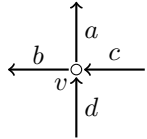
 \begin{center}
 \begin{tikzpicture}
 \node  at    (0.8,0.8) {$v$} ;

\node  at    (0.5,1.2) {$b$} ;
\node  at    (1.2,0.5) {$d$}  ;
\node  at    (1.2,1.5) {$a$}  ;
\node  at    (1.5,1.2) {$c$}  ;
\draw (1.0,1) circle [radius=0.07];

\draw[->][thick,black] (1.0,0.1) -- (1.0,0.9);
\draw[->][thick,black] (1.0,1.1) -- (1.0,1.9);
\draw[->][thick,black] (1.9,1) -- (1.1,1);
\draw[->][thick,black] (0.9,1) -- (0.1,1);
\end{tikzpicture}
\end{center}

Define $$\eta_{\!\text{\scalebox{.6}{ NE}}}(v)=\mathrm{sgn}(a,c)\, , \quad \eta_{\!\text{\scalebox{.6}{ SW}}}(v)=\mathrm{sgn}(b,d)\, ,\quad \eta_{\!\text{\scalebox{.6}{ SE}}}(v)=\mathrm{sgn}(c,d)\, ,\quad \eta_{\!\text{\scalebox{.6}{ NW}}}(v)=\mathrm{sgn}(a,b)$$
If $b$ doesn't exist (i.e. $v$ lies on the left wall) we define $\eta_{\!\text{\scalebox{.6}{ NW}}}(v)=0$ and $\eta_{\!\text{\scalebox{.6}{ SW}}}(v)=0$. If $d$ doesn't exist (i.e. $v$ lies on the bottom wall) we define $\eta_{\!\text{\scalebox{.6}{ SE}}}(v)=0$ and $\eta_{\!\text{\scalebox{.6}{ SW}}}(v)=0$.

\begin{lemma}\label{lemma1}
Suppose we have a box
\begin{center}
 \begin{tikzpicture}

\node  at    (0.3,1) {$d$} ;
\node  at    (1.7,1) {$c$}  ;
\node  at    (1,1.7) {$a$}  ;
\node  at    (1,0.7) {$b$}  ;

\node  at    (0.3,1.7) {$v$} ;
\node  at    (0.3,0.3) {$x$}  ;
\node  at    (1.7,1.7) {$w$}  ;
\node  at    (1.7,0.3) {$y$}  ;

\draw (0.5,0.5) circle [radius=0.07];
\draw (0.5,1.5) circle [radius=0.07];
\draw (1.5,0.5) circle [radius=0.07];
\draw (1.5,1.5) circle [radius=0.07];

\draw[->][thick,black] (0.5,0.6) -- (0.5,1.4);
\draw[->][thick,black] (1.5,0.6) -- (1.5,1.4);
\draw[->][thick,black] (1.4,0.5) -- (0.6,0.5);
\draw[->][thick,black] (1.4,1.5) -- (0.6,1.5);
\end{tikzpicture}
\end{center}
then $\eta_{\!\text{\scalebox{.6}{ SE}}}(y)+\eta_{\!\text{\scalebox{.6}{ NW}}}(v) \geq \eta_{\!\text{\scalebox{.6}{ NE}}}(x)+\eta_{\!\text{\scalebox{.6}{ SW}}}(w)$. \end{lemma}
\noindent\textit{Proof of Lemma \ref{lemma1}}: This can be seen by checking the different cases.\\
\text{\tiny $\bullet$}  If $\eta_{\!\text{\scalebox{.6}{ SE}}}(y)=-1$ and $\eta_{\!\text{\scalebox{.6}{ NW}}}(v)=-1$, then wlog $\tau_{b}>0$ and $\tau_{c}<0$. If $\eta_{\!\text{\scalebox{.6}{ SW}}}(w)=0$ then $\tau_{a}=0$ so $\eta_{\!\text{\scalebox{.6}{ NW}}}(v)=0$. Contradiction. If $\eta_{\!\text{\scalebox{.6}{ SW}}}(w)=1$ then $\tau_{a}<0$ and $\tau_{d}>0$ whence $\rho_{b}+\rho_{d}<\sigma_{b}+\sigma_{d}=\sigma_{a}+\sigma_{c}<\rho_{a}+\rho_{c}$. Contradiction. Therefore $\eta_{\!\text{\scalebox{.6}{ SW}}}(w)=-1\implies\tau_{a}>0\implies\tau_{d}<0\implies \eta_{\!\text{\scalebox{.6}{ NE}}}(x)=-1$, so $\eta_{\!\text{\scalebox{.6}{ SE}}}(y)+\eta_{\!\text{\scalebox{.6}{ NW}}}(v) \geq \eta_{\!\text{\scalebox{.6}{ NE}}}(x)+\eta_{\!\text{\scalebox{.6}{ SW}}}(w)$.\\
\text{\tiny $\bullet$}  If $\eta_{\!\text{\scalebox{.6}{ SE}}}(y)=-1 \mbox{ and } \eta_{\!\text{\scalebox{.6}{ NW}}}(v)=0$ then wlog $\tau_{b}>0 \mbox{ and } \tau_{c}<0 \mbox{ and } \tau_{a}=0$, whence $\eta_{\!\text{\scalebox{.6}{ SW}}}(w)=0 \mbox{ and }\eta_{\!\text{\scalebox{.6}{ NE}}}(x)=-1$, as $\eta_{\!\text{\scalebox{.6}{ NE}}}(x) \geq0 \implies \rho_{b}+\rho_{d}<\sigma_{b}+\sigma_{d}=\sigma_{a}+\sigma_{c}<\rho_{a}+\rho_{c}$. Contradiction.\\
\text{\tiny $\bullet$}  If $\eta_{\!\text{\scalebox{.6}{ SE}}}(y)=-1 \mbox{ and } \eta_{\!\text{\scalebox{.6}{ NW}}}(v)=1$ then wlog $\tau_{b}>0 \mbox{ and } \tau_{c}<0 \mbox{ and } \tau_{a}>0\mbox{ and } \tau_{d}>0$, whence $\eta_{\!\text{\scalebox{.6}{ SW}}}(w)=-1 \mbox{ and }\eta_{\!\text{\scalebox{.6}{ NE}}}(x)=1$.\\
\text{\tiny $\bullet$}  If $\eta_{\!\text{\scalebox{.6}{ SE}}}(y)=0 \mbox{ and } \eta_{\!\text{\scalebox{.6}{ NW}}}(v)=0$ then either wlog $\tau_{b}=0 \mbox{ and } \tau_{a}=0$ whence $\eta_{\!\text{\scalebox{.6}{ SW}}}(w)=0 \mbox{ and }\eta_{\!\text{\scalebox{.6}{ NE}}}(x)=0$, or wlog $\tau_{b}=0 \mbox{ and } \tau_{d}=0$ whence $\eta_{\!\text{\scalebox{.6}{ NE}}}(x)=0 \mbox{ and }\eta_{\!\text{\scalebox{.6}{ SW}}}(w)\leq0$, as $\eta_{\!\text{\scalebox{.6}{ SW}}}(w) =1 \implies$ wlog $\tau_{a}<0\mbox{ and } \tau_{c}<0$ giving $\rho_{b}+\rho_{d}=\sigma_{b}+\sigma_{d}=\sigma_{a}+\sigma_{c}<\rho_{a}+\rho_{c}$. Contradiction.\\
\text{\tiny $\bullet$}  If $\eta_{\!\text{\scalebox{.6}{ SE}}}(y)=0 \mbox{ and } \eta_{\!\text{\scalebox{.6}{ NW}}}(v)=1$ then wlog $\eta_{\!\text{\scalebox{.6}{ SW}}}(w)=0$, and we're done.\\
\text{\tiny $\bullet$} If $\eta_{\!\text{\scalebox{.6}{ SE}}}(y)=1 \mbox{ and } \eta_{\!\text{\scalebox{.6}{ NW}}}(v)=1$ then we're done too. $\square$

\begin{lemma}\label{lemma2}
Suppose we have a vertex $v \in \mathcal{V}_{\bullet}$
 \begin{center}
 \begin{tikzpicture}
 \node  at    (0.8,0.8) {$v$} ;

\node  at    (0.5,1.2) {$b$} ;
\node  at    (1.2,0.5) {$d$}  ;
\node  at    (1.2,1.5) {$a$}  ;
\node  at    (1.5,1.2) {$c$}  ;
\draw [fill] (1.0,1) circle [radius=0.07];

\draw[->][thick,black] (1.0,0.1) -- (1.0,0.9);
\draw[->][thick,black] (1.0,1.1) -- (1.0,1.9);
\draw[->][thick,black] (1.9,1) -- (1.1,1);
\draw[->][thick,black] (0.9,1) -- (0.1,1);
\end{tikzpicture}
\end{center}
then $\eta_{\!\text{\scalebox{.6}{ NE}}}(v)+\eta_{\!\text{\scalebox{.6}{ SW}}}(v) \geq \eta_{\!\text{\scalebox{.6}{ SE}}}(v)+\eta_{\!\text{\scalebox{.6}{ NW}}}(v)$. \end{lemma}
\noindent\textit{Proof of Lemma \ref{lemma2}}: Here also, we just need to check the different cases.\\
\text{\tiny $\bullet$} If $\eta_{\!\text{\scalebox{.6}{ NE}}}(v)=-1 \mbox{ and } \eta_{\!\text{\scalebox{.6}{ SW}}}(v)=-1$ then wlog $\tau_{a}>0 \mbox{ and } \tau_{c}<0$. Then $\eta_{\!\text{\scalebox{.6}{ NW}}}(v)= 0\implies \tau_{b}=0 \implies \eta_{\!\text{\scalebox{.6}{ SW}}}(v)=0$ (Contradiciton) and $\eta_{\!\text{\scalebox{.6}{ NW}}}(v)=1\implies \tau_{b}>0 \mbox{ and } \tau_{d}<0$ whence $\min(\rho_{c},\rho_{d})<\min(\sigma_{c},\sigma_{d})=\min(\sigma_{a},\sigma_{b})<\min(\rho_{a},\rho_{b})$ (Contradiction). So $\eta_{\!\text{\scalebox{.6}{ NW}}}(v)=-1$, which forces $\eta_{\!\text{\scalebox{.6}{ SE}}}(v)=-1$.\\
\text{\tiny $\bullet$} If $\eta_{\!\text{\scalebox{.6}{ NE}}}(v)=-1 \mbox{ and } \eta_{\!\text{\scalebox{.6}{ SW}}}(v)=0$ then wlog $\tau_{a}=0 \mbox{ and } \tau_{b}<0\mbox{ and } \tau_{d}>0$. Then $\eta_{\!\text{\scalebox{.6}{ NW}}}(v)= 0$, and $\eta_{\!\text{\scalebox{.6}{ SE}}}(v)\geq 0\implies \tau_{c}\geq0 \implies \min(\rho_{c},\rho_{d})<\min(\sigma_{c},\sigma_{d})=\min(\sigma_{a},\sigma_{b})<\min(\rho_{a},\rho_{b})$ (Contradiction). So $\eta_{\!\text{\scalebox{.6}{ SE}}}(v)=-1$.\\
\text{\tiny $\bullet$} If $\eta_{\!\text{\scalebox{.6}{ NE}}}(v)=0 \mbox{ and } \eta_{\!\text{\scalebox{.6}{ SW}}}(v)=0$ then either wlog $\tau_{c}=0 \mbox{ and } \tau_{b}=0$ giving $\eta_{\!\text{\scalebox{.6}{ SE}}}(v)=\eta_{\!\text{\scalebox{.6}{ NW}}}(v)=0$ or wlog  $\tau_{c}=0 \mbox{ and } \tau_{d}=0$ giving $\eta_{\!\text{\scalebox{.6}{ SE}}}(v)=0$. If $\eta_{\!\text{\scalebox{.6}{ NW}}}(v)=1$ then wlog $\tau_{a}<0 \mbox{ and } \tau_{b}<0$ whence $\min(\rho_{c},\rho_{d})=\min(\sigma_{c},\sigma_{d})=\min(\sigma_{a},\sigma_{b})<\min(\rho_{a},\rho_{b})$ (Contradiction). So $\eta_{\!\text{\scalebox{.6}{ SE}}}(v)\leq 0$.\\
\text{\tiny $\bullet$} If $\eta_{\!\text{\scalebox{.6}{ NE}}}(v)=0 \mbox{ and } \eta_{\!\text{\scalebox{.6}{ SW}}}(v)=1$ then wlog $\eta_{\!\text{\scalebox{.6}{ SE}}}(v)=0$, and we're done.\\
\text{\tiny $\bullet$} If $\eta_{\!\text{\scalebox{.6}{ NE}}}(v)=1 \mbox{ and } \eta_{\!\text{\scalebox{.6}{ SW}}}(v)=1$ then we're done.\\
Note that this inequality also holds when $v$ is on the left or bottom wall (using the definitions for $\eta_{\!\text{\scalebox{.6}{ NW}}},\eta_{\!\text{\scalebox{.6}{ SW}}},\eta_{\!\text{\scalebox{.6}{ SE}}}$ in these cases). $\square$\\

Now for each diagonal $\mathcal{D}_{i}$ we define the following sums
$$\phi_{i}:=\sum\limits_{v \in \mathcal{D}_{i}}\eta_{\!\text{\scalebox{.6}{ SW}}}(v) \quad\text{and}\quad \psi_{i}:=\sum\limits_{v \in \mathcal{D}_{i}}\eta_{\!\text{\scalebox{.6}{ NE}}}(v)\quad \text{for}\quad i=1,...,n$$ $$\beta_{i}:=\phi_{i}-\psi_{i+1}\quad \text{for}\quad i=1,...,n\quad \quad (\text{Set}\,\,\mathcal{D}_{n+1}=\emptyset\,\,\text{and}\,\,\psi_{n+1}=0)$$ Then we claim that $\beta_{i} \leq \beta_{i+1}$ for $i=1,..,n-1$. This is equivalent to showing $\phi_{i}+\psi_{i+2} \leq \phi_{i+1}+\psi_{i+1} $ for $i=1,..,n-1$, which can be seen as follows:
\begin{align*}\phi_{i}+\psi_{i+2} &= \sum\limits_{v \in \mathcal{D}_{i}}\eta_{\!\text{\scalebox{.6}{ SW}}}(v) +\sum\limits_{v \in \mathcal{D}_{i+2}}\eta_{\!\text{\scalebox{.6}{ NE}}}(v) \\
&\leq \eta_{\!\text{\scalebox{.6}{ SW}}}(v_{i\,1})+\eta_{\!\text{\scalebox{.6}{ SW}}}(v_{n\,n+1-i})+ \sum\limits_{v \in \mathcal{D}_{i+1}}\eta_{\!\text{\scalebox{.6}{ SE}}}(v)+\eta_{\!\text{\scalebox{.6}{ NW}}}(v)\quad-\eta_{\!\text{\scalebox{.6}{ SE}}}(v_{n\,n-i})-\eta_{\!\text{\scalebox{.6}{ NW}}}(v_{i+1\,1}) \quad \text{(Lemma \ref{lemma1})} \\
&=\sum\limits_{v \in \mathcal{D}_{i+1}}\eta_{\!\text{\scalebox{.6}{ SE}}}(v)+\eta_{\!\text{\scalebox{.6}{ NW}}}(v)\\
&\leq\sum\limits_{v \in \mathcal{D}_{i+1}}\eta_{\!\text{\scalebox{.6}{ NE}}}(v)+\eta_{\!\text{\scalebox{.6}{ SW}}}(v)\quad \text{(Lemma \ref{lemma2})} \\
&=\phi_{i+1}+\psi_{i+1}
\end{align*}
Note that $\beta_{n}=\eta_{\!\text{\scalebox{.6}{ SW}}}(v_{n\,n})=0$, so we therefore have $0=\beta_{n}\geq \beta_{n-1} \geq ...\geq \beta_{1}$.\\

We now make the following claim: $$\beta_{1} \geq 0\text{ with equality iff }\tau_{a}=0\text{ for all arrows }a \text{ that touch a star vertex}$$ This claim would then complete the proof as it forces $\beta_{1} = 0$ and thus implies $T_{v}=S_{v}$ for all sub-diagonal vertices. So we would get two tropical critical points for a diagram of size $n-1$, which must be equal be the induction hypothesis.

To prove the claim first note that if $v$ is a vertex on the sub-diagonal $\mathcal{D}_{2}$ 
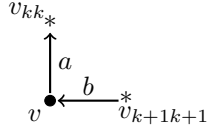
\begin{center}
 \begin{tikzpicture}
 \node  at    (0.8,0.8) {$v$} ;
 \node  at    (0.7,2.2) {$v_{kk}$} ;
 \node  at    (2.5,0.8) {$v_{k+1k+1}$} ;

\node  at    (1.2,1.5) {$a$}  ;
\node  at    (1.5,1.2) {$b$}  ;
\draw [fill] (1.0,1) circle [radius=0.07];

\node  at    (1.0,2)  {*}  ;
\node  at    (2,1)  {*}  ;
\draw[->][thick,black] (1.0,1.1) -- (1.0,1.9);
\draw[->][thick,black] (1.9,1) -- (1.1,1);
\end{tikzpicture}
\end{center}
then we have $\eta_{\!\text{\scalebox{.6}{ NE}}}(v)\geq0  \implies \tau_{a}=\tau_{b}=0$ since otherwise we would have wlog $\sigma_{a}\geq\rho_{a} \mbox{ and } \sigma_{b}>\rho_{b}$ which implies the contradiction $\lambda_{k}-\lambda_{k+1}=\sigma_{a}+\sigma_{b}>\rho_{a}+\rho_{b}=\lambda_{k}-\lambda_{k+1}$. Hence we have $\eta_{\!\text{\scalebox{.6}{ NE}}}(v)\leq0 $ for all $v \in \mathcal{D}_{2}$.

Now let $H$ be the subset of the star vertices for which $\eta_{\!\text{\scalebox{.6}{ SW}}}(v)=-1$. Suppose $H=\emptyset $. Then $\beta_{1} \geq 0$ and $\beta_{1} = 0 \implies \eta_{\!\text{\scalebox{.6}{ NE}}}(v)=0 \quad \forall  v \in\mathcal{D}_{2} $. Hence $\beta_{1} = 0$ implies $\tau_{a}=0$ for all arrows $a$ that touch a star vertex. Suppose $H\neq\emptyset$. Let $H'$ be the subset of sub-diagonal vertices that are connected to a star vertex in $H$ by an arrow. Then $\eta_{\!\text{\scalebox{.6}{ NE}}}(v)=-1$ for all $v \in H'$ since if we have $\eta_{\!\text{\scalebox{.6}{ NE}}}(v)\geq0$ then there is an arrow $a$ joining it to a star vertex $w\in H$ such that $\tau_{a}=0$, implying the contradiction $\eta_{\!\text{\scalebox{.6}{ SW}}}(w)=0$. Also note that $\vert H' \vert > \vert H \vert$. Hence $$\beta_{1}=\sum\limits_{v \in \mathcal{D}_{1}}\eta_{\!\text{\scalebox{.6}{ SW}}}(v) -\sum\limits_{v \in \mathcal{D}_{2}}\eta_{\!\text{\scalebox{.6}{ NE}}}(v)\geq \sum\limits_{v \in H}\eta_{\!\text{\scalebox{.6}{ SW}}}(v)-\sum\limits_{v\in H'}\eta_{\!\text{\scalebox{.6}{ NE}}}(v)=\vert H' \vert - \vert H \vert >0$$
If $\beta_{1} = 0$ then we must have $\vert H' \vert = \vert H \vert$, so $H=\emptyset$ and again $\tau_{a}=0$ for all arrows $a$ that touch a star vertex, hence the claim is true. This completes the proof of Proposition 5.6, and so of uniqueness.\\

\subsection*{The point $p_{\lambda}$ lies in $Z_{t^{\lambda}}(\mathcal{K}_{>0})^+$} Take $\tilde{\lambda}$ a lift of $\lambda$. Suppose in the special chart the critical point is given by $(x_v)_{v\in\mathcal{V}}\in (\mathcal{K}_{>0})^{\mathcal{V}}$. Let $z_a=\frac{x_{\mathrm{h}(a)}}{x_{\mathrm{t}(a)}}$ and let $\delta_v=\mathrm{val}(x_v)$ and $\sigma_a=\mathrm{val}(z_a)$. Then we need to show that $\min\limits_{a \in \mathcal{A}}(\sigma_{a}) \geq 0$. Suppose $\sigma_{a}<0$ for some $a \in \mathcal{A}$, then since the $\sigma_{a}$ obey the tropical critical point conditions $$\min\limits_{a:h(a)=v}\sigma_{a}=\min\limits_{a:t(a)=v}\sigma_{a}\quad \forall v \in \mathcal{V}^{\bullet}$$ we can find a path of arrows $\pi$ between two star vertices $v_{s}$ and $v_{e}$ such that $\sigma_{a}<0$ for all $a \in \pi$. Then $\delta_{v_s}>\delta_{v_e}$, which contradicts the fact that $\tilde{\lambda}$ was dominant.\\

\subsection*{The point $p_{\lambda}$ has ``weight zero''} Next we prove that $\mathrm{wt}(p_{\lambda})=1\in T^{\vee}(\mathcal{K})$. Again take $\tilde{\lambda}$ a lift of $\lambda$ and suppose in the special chart the critical point is given by $(x_v)_{v\in\mathcal{V}}\in (\mathcal{K}_{>0})^{\mathcal{V}}$. Let $z_a=\frac{x_{\mathrm{h}(a)}}{x_{\mathrm{t}(a)}}$ and let $\delta_v=\mathrm{val}(x_v)$ and $\sigma_a=\mathrm{val}(z_a)$. In the special chart, $\mathrm{wt}$ is given by $\tilde{\gamma}$. Hence we need to show $\frac{\zeta_{i}}{\zeta_{i+1}}=\frac{\zeta_{i-1}}{\zeta_{i}}$ for all $2 \leq i \leq n$. We will use the following lemma:

\begin{lemma}\label{lemmawtzero}
Suppose we have a quiver
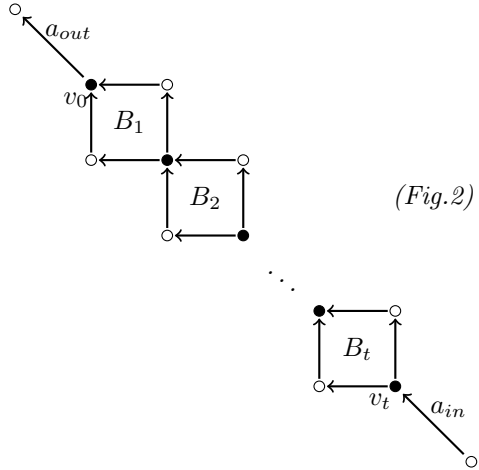
\begin{center}
 \begin{tikzpicture}
 \node  at    (5.5,3.5)  {(Fig.2)};

 \node[fill=white] at    (1.5,4.5) {$B_{1}$}  ;
\node[fill=white] at    (2.5,3.5) {$B_{2}$}  ;
\node[fill=white] at    (4.5,1.5) {$B_{t}$}  ;

\node[fill=white] at    (0.7,5.7) {$a_{out}$}  ;
\node[fill=white] at    (5.7,0.7) {$a_{in}$}  ;
\node[fill=white] at    (0.8,4.8) {$v_{0}$}  ;
\node[fill=white] at    (4.8,0.8) {$v_{t}$}  ;
\node[fill=white] at    (3.5,2.5) {$\ddots$}  ;
\draw        (0.0,6) circle [radius=0.07];
\draw [fill] (1.0,5) circle [radius=0.07];
\draw  (1.0,4) circle [radius=0.07];
\draw        (2.0,5) circle [radius=0.07];
\draw [fill] (2.0,4) circle [radius=0.07];
\draw  (2.0,3) circle [radius=0.07];
\draw        (3.0,4) circle [radius=0.07];
\draw     [fill]   (3.0,3) circle [radius=0.07];
\draw     [fill]   (4.0,2) circle [radius=0.07];
\draw        (4.0,1) circle [radius=0.07];
\draw        (5.0,2) circle [radius=0.07];
\draw      [fill]  (5.0,1) circle [radius=0.07];
\draw        (6.0,0) circle [radius=0.07];

\draw[->][thick,black] (5.9,0.1) -- (5.1,0.9);
\draw[->][thick,black] (1,4.1) -- (1,4.9);
\draw[->][thick,black] (2,4.1) -- (2,4.9);
\draw[->][thick,black] (2,3.1) -- (2,3.9);
\draw[->][thick,black] (3,3.1) -- (3,3.9);
\draw[->][thick,black] (4,1.1) -- (4,1.9);
\draw[->][thick,black] (5,1.1) -- (5,1.9);
\draw[->][thick,black] (1.9,5) -- (1.1,5);
\draw[->][thick,black] (1.9,4) -- (1.1,4);
\draw[->][thick,black] (2.9,4) -- (2.1,4);
\draw[->][thick,black] (2.9,3) -- (2.1,3);
\draw[->][thick,black] (4.9,2) -- (4.1,2);
\draw[->][thick,black] (4.9,1) -- (4.1,1);
\draw[->][thick,black] (0.9,5.1) -- (0.1,5.9);
\end{tikzpicture}
\end{center}
with a variable $z_{a}$ attached to each arrow $a$ such that the \textit{tropical box relations} $z_{a_{1}}z_{a_{2}}=z_{a_{3}}z_{a_{4}}$ hold whenever $a_{1},a_{2},a_{3},a_{4}$ form a square
and the critical conditions hold at each black vertex. For $1 \leq j \leq t$ let $O_{j}=z_{a_{1}}z_{a_{3}}$ where the box $B_{j}$ is given by \begin{center}
 \begin{tikzpicture}
\node  at    (3.5,1)  {(Fig.3)};
\node  at    (0.3,1) {$a_{2}$} ;
\node  at    (1.7,1) {$a_{3}$}  ;
\node  at    (1,1.7) {$a_{4}$}  ;
\node  at    (1,0.7) {$a_{1}$}  ;
\draw[->][thick,black] (0.5,0.6) -- (0.5,1.4);
\draw[->][thick,black] (1.5,0.6) -- (1.5,1.4);
\draw[->][thick,black] (1.4,0.5) -- (0.6,0.5);
\draw[->][thick,black] (1.4,1.5) -- (0.6,1.5);
\end{tikzpicture}
\end{center}
Similarly, let $I_{j}=z_{a_{2}}z_{a_{4}}$.
Let $K_{t}=\prod z_{a}$ where the product is over a (any) path from $v_{t}$ to $v_{0}$. Then we have
$$\prod^{t}_{j=1} O_{j}.\frac{z_{a_{out}}}{z_{a_{in}}}=K_{t}\qquad\text{ and }\qquad\prod^{t}_{j=1} I_{j}.\frac{z_{a_{out}}}{z_{a_{in}}}=K_{t}$$ Note this agrees with $\prod\limits^{t}_{j=1} O_{j}\prod\limits^{t}_{j=1} I_{j}=K_{t}^2$.
\end{lemma}
\noindent\textit{Proof of Lemma \ref{lemmawtzero}}: We prove this by induction on $t$. For $t=0$ it is clear so assume $t>0$. Suppose $B_{t}$ looks like \textit{Fig.3} and let $x=z_{a_{1}}z_{a_{2}}=z_{a_{3}}z_{a_{4}}$. Then
\begin{align*}
\prod^{t}_{j=1}O_{j}.\frac{z_{a_{out}}}{z_{a_{in}}} & =\prod^{t-1}_{j=1}O_{j}.z_{a_{out}}\frac{z_{a_{1}}z_{a_{3}}}{z_{a_{1}}+z_{a_{3}}} \quad\quad \text{since }z_{a_{in}}=z_{a_{1}}+z_{a_{3}}\\
& =\prod^{t-1}_{j=1}O_{j}.z_{a_{out}}\frac{x}{z_{a_{2}}+z_{a_{4}}} \quad\quad\quad \text{by the box relation }\\
& =K_{t-1}.x \qquad \text{by induction hypothesis}\\
& =K_{t}
\end{align*}
\noindent(The statement for $I_{j}$ also follows from this). $\square$\\

Now consider all the arrows with either head or tail in the $i^{th}$ diagonal $\mathcal{D}_{i}$. These form a diagram like in \textit{Fig.2} and $\frac{\zeta_{i+1}\zeta_{i-1}}{\zeta_{i}^{2}}=\frac{\prod_{j} O_{j}.\frac{z_{a_{out}}}{z_{a_{in}}}}{K}$ which equals 1 by Lemma \ref{lemmawtzero} \\
\begin{example}
Let $n=3$ and $\tilde{\lambda}=(3,1,0)$. Then the critical point is given in the special coordinates by:
\end{example}
\begin{center}
 \begin{tikzpicture}[scale=2]
\draw [fill] (0.5,0.5) circle [radius=0.04];
\draw [fill] (0.5,1.5) circle [radius=0.04];
\draw [fill] (1.5,0.5) circle [radius=0.04];
\node  at    (0.5,2.5) {*} ;
\node  at    (1.5,1.5) {*}  ;
\node  at    (2.5,0.5) {*}  ;

\node[fill=white] at (0.1,1) {\scalebox{0.6}{$t^{5/6}-\frac{t^{7/6}}{2}+\frac{3t^{9/6}}{8}-\frac{5t^{11/6}}{16}+\frac{35t^{13/6}}{128}\dots$}}  ;
\node[fill=white] at (0.1,2) {\scalebox{0.6}{$t^{5/6}+\frac{t^{7/6}}{2}-\frac{t^{9/6}}{8}+\frac{t^{11/6}}{16}-\frac{t^{13/6}}{128}\dots$}}  ;
\node[fill=white] at (2.4,1) {\scalebox{0.6}{$t^{3/6}-\frac{t^{5/6}}{2}+\frac{3t^{7/6}}{8}-\frac{5t^{9/6}}{16}+\frac{35t^{11/6}}{128}\dots$}}  ;
\node[fill=white] at (0.5,0.3) {\scalebox{0.6}{$t^{5/6}-\frac{t^{7/6}}{2}+\frac{3t^{9/6}}{8}-\frac{5t^{11/6}}{16}+\frac{35t^{13/6}}{128}\dots$}}  ;
\node[fill=white] at (1.1,1.7) {\scalebox{0.6}{$t^{7/6}-\frac{t^{9/6}}{2}+\frac{3t^{11/6}}{8}\dots$} } ;
\node[fill=white] at (2.4,0.3) {\scalebox{0.6}{$t^{3/6}+\frac{t^{5/6}}{2}-\frac{t^{7/6}}{8}+\frac{t^{9/6}}{16}-\frac{5t^{11/6}}{128}\dots$}}  ;

\node[fill=white] at    (0.7,2.6) {$t^3$}  ;
\node[fill=white] at    (1.7,1.6) {$t^1$}  ;
\node[fill=white] at    (2.7,0.6) {$t^0$}  ;

\draw[->][thick,black] (0.5,0.6) -- (0.5,1.4);
\draw[->][thick,black] (0.5,1.6) -- (0.5,2.4);
\draw[->][thick,black] (1.5,0.6) -- (1.5,1.4);
\draw[->][thick,black] (2.4,0.5) -- (1.6,0.5);
\draw[->][thick,black] (1.4,0.5) -- (0.6,0.5);
\draw[->][thick,black] (1.4,1.5) -- (0.6,1.5);
\end{tikzpicture}\\
\end{center}
\bigbreak

\section{Integrality of $p_{\lambda}$}
Next we want to study the integrality of the critical point $p_{\lambda}$.
\begin{defn}
Given $\lambda\in P^+$, we say $p_{\lambda}$ is integral if it lies in $Z(\mathbb{C}((t)))$. Let $\mathcal{P}$ be the subset of $\lambda\in P^+$ such that $p_{\lambda}$ is integral.
\end{defn}
First note that $p_{\lambda}$ is integral iff its representation in one of the toric charts $\theta\in\Theta$ is integral. Suppose $z_{\lambda}\in(\mathcal{K}_{>0})^{\mathcal{A}}$ is the representation of $p_{\lambda}$ in the chart $\theta_{\mathcal{M}}$, and let $\sigma_{\lambda}\in(\mathbb{R})^{\mathcal{A}}$ be its valuation. Then it follows from the construction of the critical point in section 5 that $z_{\lambda}$ is integral iff $\sigma_{\lambda}$ is integral. Recall $\sigma_{\lambda}$ is the unique solution to the tropical critical conditions (5.2) with highest weight $\lambda$.

Now we introduce a combinatorial object which describes the solutions to the tropical critical conditions. Let $\mathcal{U}$ be the grid of $n(n-1)/2$ boxes in upper triangular form. Consider \textit{fillings} of $\mathcal{U}$ where we put a non-negative real number in each box. For example if $n=4$ it looks like:
\begin{center}
 \begin{tikzpicture}[scale=0.75]
\node[fill=white] at    (0.5,2.5) {$n_{12}$}  ;
\node[fill=white] at    (1.5,2.5) {$n_{13}$}  ;
\node[fill=white] at    (1.5,1.5) {$n_{23}$}  ;
\node[fill=white] at    (2.5,2.5) {$n_{14}$}  ;
\node[fill=white] at    (2.5,1.5) {$n_{24}$}  ;
\node[fill=white] at    (2.5,0.5) {$n_{34}$}  ;
\draw[thick,black] (0,2) -- (0,3) -- (3,3) -- (3,0) -- (2,0) -- (2,3)--(1,3)--(1,1)--(3,1)--(3,2)--(0,2);
 \end{tikzpicture}
\end{center}
We say a filling $\{n_{ij}\}_{1 \leq i <j \leq n} $ is \textit{ideal } if $n_{ij}=\max\{n_{i+1\,j},n_{i\,j-1}\}$ for $j-i \geq 2$ and \textit{integral} if all the $n_{ij}$ are integral. Note that an ideal filling is determined by its values on the first diagonal since $n_{ij}=\max\limits_{i \leq k \leq j-1}\{n_{kk+1}\}$.

Let $\lambda \in P^+$. We say $\{n_{ij}\}_{1 \leq i <j \leq n} $ is an ideal fillings for $\lambda$ if it is an ideal filling and $ \sum n_{ij} \alpha_{ij}=\lambda$.
Then we have the following proposition.

\begin{prop}
Let $ \lambda \in P^+$. Then we have a bijective correspondence between $$ \{\text{solutions to the tropical critical conditions with highest weight } \lambda \} \leftrightarrow \{\text{ideal fillings for }\lambda\}.$$ Suppose $(\sigma_{a})_{a\in \mathcal{A}}$ is a solution to the tropical critical conditions with highest weight $\lambda$. Recall the map $\pi:\mathcal{V}_{\bullet}\rightarrow \mathbb{Q}$ given by $\pi(v):= \min\limits_{a:\rm{h}(a)=v}\sigma_{a}= \min\limits_{a:\rm{t}(a)=v}\sigma_{a}$. Then the bijective correspondence takes this solution to the tropical critical conditions to the ideal filling for $\lambda$ given by setting $n_{ij}=\pi(v_{ji})$ for $1 \leq i < j \leq n$. Furthermore, the bijection preserves integrality.
\end{prop}
\begin{remark}
Of course, we know from the previous section that both these sets just consist of a single point.
\end{remark}
\begin{example}
Suppose $n=3$ and $\lambda$ is given by $(6,3,-2)$, then we can compute the unique solution to the tropical critical conditions and the unique ideal filling to be as follows. Note that $\lambda=\frac{3}{2}\alpha_{12}+\frac{13}{6}\alpha_{23}+\frac{13}{6}\alpha_{13}$.
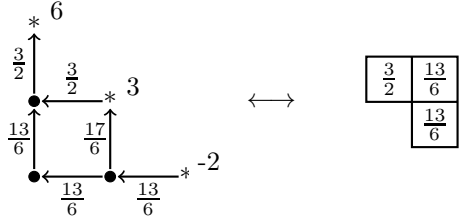
\begin{center}
 \begin{tikzpicture}
\draw [fill] (0.5,0.5) circle [radius=0.07];
\draw [fill] (0.5,1.5) circle [radius=0.07];
\draw [fill] (1.5,0.5) circle [radius=0.07];
\node  at    (0.5,2.5) {*} ;
\node  at    (1.5,1.5) {*}  ;
\node  at    (2.5,0.5) {*}  ;
\node[fill=white] at    (0.3,1) {$\frac{13}{6}$}  ;
\node[fill=white] at    (0.3,2) {$\frac{3}{2}$}  ;
\node[fill=white] at    (1.3,1) {$\frac{17}{6}$}  ;
\node[fill=white] at    (1,0.2) {$\frac{13}{6}$}  ;
\node[fill=white] at    (1,1.8) {$\frac{3}{2}$}  ;
\node[fill=white] at    (2,0.2) {$\frac{13}{6}$}  ;
\node[fill=white] at    (0.8,2.7) {6}  ;
\node[fill=white] at    (1.8,1.7) {3}  ;
\node[fill=white] at    (2.8,0.7) {-2}  ;
\draw[->][thick,black] (0.5,0.6) -- (0.5,1.4);
\draw[->][thick,black] (0.5,1.6) -- (0.5,2.4);
\draw[->][thick,black] (1.5,0.6) -- (1.5,1.4);
\draw[->][thick,black] (2.4,0.5) -- (1.6,0.5);
\draw[->][thick,black] (1.4,0.5) -- (0.6,0.5);
\draw[->][thick,black] (1.4,1.5) -- (0.6,1.5);
\end{tikzpicture}
\begin{tikzpicture}
\node  at    (1.5,2.5) {$\longleftrightarrow $} ;
\node  at    (2.5,1) {} ;
\end{tikzpicture}
\begin{tikzpicture}[scale=0.6]
\draw[thick,black] (2,3) -- (0,3) -- (0,4) -- (2,4) -- (2,2) -- (1,2) -- (1,4);
\node  at    (0.5,3.5) {$\frac{3}{2}$} ;
\node  at    (1.5,2.5) {$\frac{13}{6}$} ;
\node  at    (1.5,3.5) {$\frac{13}{6}$} ;
\node  at    (1.5,0.5) {} ;
 \end{tikzpicture}
 \end{center}
\end{example}
\begin{remark}
This theorem then gives a simpler way to determine whether the critical point for $\lambda$ is integral or not.
\end{remark}

The uniqueness of the ideal filling for $\lambda$ also gives us the following corollary.
\begin{cor}
Given $\lambda\in P^+$, there exists a unique way to write $$\lambda=\sum_{P\subset G}c_{P}\lambda_{P}$$ such that $P\subset G$ is a parabolic subgroup, $c_{P}\in \mathbb{R}_{\geq0}$ and the set of parabolics $\{P:c_{P}\neq 0\}$ form a chain w.r.t inclusion. Furthermore, $p_{\lambda}$ is integral iff all the coefficients $c_{P}$ are.
\end{cor}
\noindent\textit{Proof of Corollary 6.6}:
Suppose the unique ideal filling for $\lambda$ is $\{n_{ij}\}_{1 \leq i <j \leq n} $. Define $I_{1}=I$ and define $I_{k}$ inductively as follows: Let $m_{k}=\min\limits_{i\in I_{k}}n_{i i+1}$ and let $I_{k+1}=\{i\in I_k:n_{i i+1}>m_k\}$. This algorithm stops when $I_{k+1}= \emptyset$. Let $P_k$ be the parabolic subgroup such that $I^{P}=I_{k}$. Then the $P_{k}$ form a chain w.r.t inclusion and $\lambda=\sum\limits_{k\geq1}(m_{k}-m_{k-1})\lambda_{P_k}$. ($m_0=0$). To see this, first observe that $\lambda_{P_{k+1}}=\sum\limits_{\alpha_{ij}:n_{ij}>m_{k}}\alpha_{ij}$ since $n_{ij}>m_k$ iff at least one of $n_{ii+1},...,n_{j-1j}$ is $>m_k$ iff the interval $[i,j-1]$ intersects $I_{k+1}$. Then
\begin{align*}
\lambda&= \sum n_{ij} \alpha_{ij}\\
&=\sum (n_{ij}-m_1) \alpha_{ij}+m_1\sum \alpha_{ij}\\
&=\sum\limits_{\alpha_{ij}:n_{ij}>m_1} (n_{ij}-m_1) \alpha_{ij}+m_1\lambda_{P_1}\\
&=\sum\limits_{\alpha_{ij}:n_{ij}>m_1} (n_{ij}-m_2) \alpha_{ij}+m_1\lambda_{P_1}+(m_2-m_1)\sum\limits_{\alpha_{ij}:n_{ij}>m_1}\alpha_{ij}\\
&=\sum\limits_{\alpha_{ij}:n_{ij}>m_1} (n_{ij}-m_2) \alpha_{ij}+m_1\lambda_{P_1}+(m_2-m_1)\lambda_{P_2}\\
& \dots\\
&=\sum\limits_{k\geq1}(m_{k}-m_{k-1})\lambda_{P_k}&&(6.1)
\end{align*}
Conversely, suppose we have $\lambda=\sum_{k=1}^{K}c_{k}\lambda_{Q_k}$ with $c_{k}\neq 0$ for $k=1,..,K$ and $Q_1\subset Q_2 \subset ...\subset Q_K$.   Given a parabolic subgroup $P$ we get a corresponding ideal filling $F_P$ for $\lambda_P$, namely the one determined by $$n_{i\,i+1}=\begin{cases} 1 \quad i \in I^P \\ 0  \quad i \notin I^P  \end{cases}$$ Hence the point-wise sum $\sum_{k=1}^{K}c_{k}F_{Q_k}$ gives an ideal filling for $\lambda$.

Let $\{n_{ij}\}_{1 \leq i <j \leq n} $ be the unique ideal filling for $\lambda$ and $\lambda=\sum_{k=1}^{L}d_{k}\lambda_{P_k}$ with $d_{k}\neq 0$ $k=1,..,L$ and $P_1\subset P_2 \subset ...\subset P_K$ the decomposition into parabolics constructed above (6.1). We need to show $P_k=Q_k$ and $c_k=d_k$ for all $k$. Suppose $i\in I^{Q_1}\backslash I^{P_1}$, then $i\in I^{Q_1}$ implies $n_{ii+1}\neq 0$ while $i\notin I^{P_1}$ implies $n_{ii+1}= 0$, so $ I^{Q_1}\subset I^{P_1}$, and similarly $I^{P_1}\subset I^{Q_1}$. To show $c_1=d_1$, wlog assume $c_1<d_1$ and let $i\in  I^{P_1}\backslash I^{P_2}$. Then according to the decomposition $\{P_k\}$ we must have $n_{ii+1}=d_1$, but according to the decomposition $\{Q_k\}$, $n_{ii+1}\geq c_1$, so we must have $c_1=d_1$. Then we can apply the same argument to $P_2$ and $Q_2$ and so on to get that $P_k=Q_k$ and $c_k=d_k$ for all $k$. $\square$\\

\noindent\textit{Proof of Proposition 6.2}:
We will define maps going each way between the two sets and show they are well defined and inverse to each other. As before we will use the coordinates given by the chart $\theta_{\mathcal{M}}$. First we define the map from $\{\textit{ideal fillings for }\lambda\}$ to $\{\textit{tropical critical points for } \lambda \} $ and show it is well defined.\\

\noindent\underline{Map from ideal fillings to solutions to the tropical critical conditions}

Let $\{n_{ij}\}_{1 \leq i <j \leq n}$ be an ideal filling for $\lambda$. Let $\tilde{\lambda}=(\lambda_{1},...,\lambda_{n})$ be the lift of $\lambda$ for which $\sum\lambda_{i}=0$. For $1 \leq i \leq j \leq n$ define $H^{h}_{ij}=\sum\limits_{k>j}n_{ik}$ and $H^{v}_{ij}=\sum\limits_{k<i}n_{kj}$. Let $\delta_{v_{ji}}=H^{h}_{ij}-H^{v}_{ij}$ for $1 \leq i \leq j \leq n$. We will show that this defines a solution to the tropical critical conditions for $\tilde\lambda$ (and hence for $\lambda$).

We can compute the corresponding arrow coordinates as follows. For $1\leq i \leq j<n$ we have \begin{align*}\delta_{v_{ji}}-\delta_{v_{j+1\,i}}&=(H^{h}_{i\,j}-H^{v}_{i\,j})-(H^{h}_{i\,j+1}-H^{v}_{i\,j+1})\\
& = (H^{h}_{i\,j}-H^{h}_{i\,j+1})-(H^{v}_{i\,j}-H^{v}_{i\,j+1})\\
& = H^{v}_{i+1\,j+1}-H^{v}_{i\,j}\end{align*}
Similarly for the horizontal arrows, if $1\leq i < j\leq n$ we have \begin{align*}\delta_{v_{ji}}-\delta_{v_{j\,i+1}}&=(H^{h}_{i\,j}-H^{v}_{i\,j})-(H^{h}_{i+1\,j}-H^{v}_{i+1\,j})\\
&= (H^{h}_{i\,j}-H^{h}_{i+1\,j})-(H^{v}_{i\,j}-H^{v}_{i+1\,j})\\
& = H^{h}_{i\,j-1}-H^{h}_{i+1\,j}\end{align*}
Note that both these expressions are $\geq 0$, so already we get that the point lies in $\{\mathcal{W}^{t}\geq0\}$.

Next recall the $\mu_{k}\in X^{*}(T_{GL_{n}^{\vee}})$ defined in Section 2. Then we have
\begin{align*}\lambda_{k}
&=\langle \tilde{\lambda},\mu_{k} \rangle\\
&=\langle \sum n_{ij}(\epsilon_{i}-\epsilon_{j}),\mu_{k}\rangle\\
&=\sum\limits_{l>k}n_{kl}-\sum\limits_{l<k}n_{lk}=H^{h}_{k\,k}-H^{v}_{k\,k}\\
&=\delta_{v_{kk}}\end{align*}
Hence the tropical point we have defined lies in the fiber over $\lambda$.

Now we show this point satisfies the tropical critical point conditions. Firstly let $\bar{H}^{v}_{i\,j}=H^{v}_{i\,j}+n_{ij}$ and $\bar{H}^{h}_{i\,j}=H^{h}_{i\,j}+n_{ij}$ for $1 \leq i < j \leq n$. Then we have
\begin{lemma}\label{ltech}
Let $j-i\geq1$, then at least one of  $\bar{H}^{v}_{i\,j}=\bar{H}^{v}_{i\,j+1}$ and $\bar{H}^{h}_{i\,j+1}=\bar{H}^{h}_{i+1\,j+1}$ must be true. Hence we have $\min\{\bar{H}^{v}_{i\,j+1}-\bar{H}^{v}_{i\,j},\bar{H}^{h}_{i\,j+1}- \bar{H}^{h}_{i+1\,j+1}\}=0$.
\end{lemma}
\noindent\textit{Proof of Lemma \ref{ltech}}: We have $\bar{H}^{v}_{i\,j}\leq\bar{H}^{v}_{i\,j+1}$ so suppose $\bar{H}^{v}_{ij}<\bar{H}^{v}_{i\,j+1}$. Then $\exists l$ with $1 \leq l \leq i$ and $n_{lj}<n_{lj+1}$. Hence $\max\limits_{l \leq k \leq j-1}\{n_{kk+1}\}<\max\limits_{l \leq k\leq j}\{n_{kk+1}\}$, and so $n_{jj+1}>\max\limits_{l \leq k\leq j-1}\{n_{kk+1}\}$. In particular, since $l \leq i \leq j-1$ we see $n_{jj+1}>n_{ii+1}$. Then for $k\geq j+1$ we have $$n_{ik}=\max\limits_{i \leq m\leq k-1}\{n_{mm+1}\}=\max\limits_{i+1 \leq m\leq k-1}\{n_{mm+1}\}=n_{i+1\,k}.$$ Hence $\bar{H}_{i\,j+1}^{h}-\bar{H}_{i+1\,j+1}^{h}=0$. This completes the proof of Lemma \ref{ltech}. $\square$\\

Now let $v_{ji}\in \mathcal{V}_{\bullet}$ with $1<i<j<n$, i.e. $v_{ji}$ doesn't lie on either wall. Then the minimum over incoming arrows to $v_{ji}$ is
\begin{equation}\label{eq1}
\begin{split}
\min\{H^{v}_{i+1\,j+1}-H^{v}_{i\,j},H^{h}_{i\,j-1}-H^{h}_{i+1\,j}\}&=n_{ij}+\min\{\bar{H}^{v}_{i\,j+1}-\bar{H}^{v}_{i\,j},\bar{H}^{h}_{i\,j+1}-\bar{H}^{h}_{i+1\,j+1}\}\\
&=n_{ij}
\end{split}
\end{equation}
and the minimum over outgoing arrows from $v_{ji}$ is
\begin{equation}\label{eq2}
\begin{split}
\min\{H^{v}_{i+1\,j}-H^{v}_{i\,j-1},H^{h}_{i-1\,j-1}-H^{h}_{i\,j}\}&=n_{ij}+\min\{\bar{H}^{v}_{i\,j+1}-\bar{H}^{v}_{i\,j},\bar{H}^{h}_{i\,j+1}-\bar{H}^{h}_{i+1\,j+1}\}\\
&=n_{ij}
\end{split}
\end{equation}
So the tropical critical point condition is satisfied at $v_{ji}$. In the case $v_{ji}$ lies on the left wall, i.e $i=1$, we have only one outgoing arrow $H^{v}_{i+1\,j}-H^{v}_{i\,j-1}=n_{1\,j}$ and if $j=n$ we have only one incoming arrow $H^{h}_{i\,n-1}-H^{h}_{i+1\,n}=n_{i\,n}$, so holds in these cases too. Thus we have defined a solution to the tropical critical conditions for $\lambda$ starting with an ideal filling for $\lambda$.

\begin{remark}
If we let $H_{ij}=H^{h}_{ij}+H^{v}_{ij}$ and $\bar{H}_{i\,j}=H_{i\,j}+n_{ij}$, then the above lemma shows that we have $\max\{H_{i\,j-1},H_{i+1\,j}\}=\bar{H}_{ij}$ for $j-i \geq 2$. This is clear since \begin{align*}\bar{H}_{i\,j}-\max\{H_{i\,j-1},H_{i+1\,j}\}&=\min\{\bar{H}_{i\,j}-H_{i\,j-1},\bar{H}_{i\,j}-H_{i+1\,j}\}\\
&=\min\{H^{v}_{i\,j}-H^{v}_{i\,j-1},H^{h}_{i\,j}-H^{h}_{i+1\,j}\}=0\end{align*}
\end{remark}

\noindent\underline{Map from solutions to the tropical critical conditions to ideal fillings}

For the inverse, suppose a solution to the tropical critical conditions for $\lambda$ is given by $(\sigma_{a})_{a\in \mathcal{A}}$ in the chart $\theta_{\mathcal{M}}$. For $v \in \mathcal{V}_{\bullet}$ let $\mathrm{inc}(v)$ be the set of arrows that are touching $v$ and let $\pi(v)= \min\limits_{a \in \rm{inc}(a)}\sigma_{a}$. Then set $n_{ij}=\pi(v_{ji})$ for $1 \leq i < j \leq n$. We will show that this defines an ideal filling for $\lambda$.

Let $\tilde{\lambda}=(\lambda_{1},...,\lambda_{n})$ be the lift of $\lambda$ such that $\sum\lambda_{i}=0$ and let $(\delta_{v})_{v \in \mathcal{V}}$ be the vertex coordinates of the solution to the tropical critical conditions for $\tilde{\lambda}$. Note that by the weight zero property of the tropical critical point we have $\delta_{v}=0$ for the bottom left vertex. The following lemma gives an expression for $\delta_{v}$ in terms of the $\sigma_{a}$.
\begin{lemma}\label{linv}
For $v \in \mathcal{V}$ we have $$\delta_{v}=\sum\limits_{w\in \mathrm{bel}(v)}\pi(w)-\sum\limits_{w \in \mathrm{lef}(v)}\pi(w)$$
where $\mathrm{bel}(v)$ is the set of vertices directly below $v$ and $\mathrm{lef}(v)$ is the set of vertices directly to the left of $v$.
\end{lemma}
\noindent\textit{Proof of Lemma \ref{linv}}:
This is true for the bottom left vertex as both sides of the equation are zero. We will prove this by inducting from the bottom left by horizontal and vertical arrows. We prove the inductive step for horizontal induction, the vertical induction is proved similarly. Suppose we have
 \begin{center}
 \begin{tikzpicture}[thick, scale=1.6]
\node  at    (0,0) {$v$} ;
\node  at    (1,0) {$w$}  ;
\node[fill=white] at    (0.5,0.13) {$c$}  ;
\draw[->][thick,black] (0.8,0) -- (0.2,0);
\end{tikzpicture}
\end{center}
and suppose the lemma holds for $v$ then we want to prove it for $w$.

Firstly suppose the diagram below is a part of the full diagram for a solution to the tropical critical conditions, where the arrow $c_{m}$ lies in the bottom wall of the full diagram.
\begin{center}
 \begin{tikzpicture}[thick, scale=1.6]
\node  at    (0,0) {$w_{m}$} ;
\node  at    (0,1) {$w_{m-1}$}  ;
\node  at    (0,2) {$w_{2}$}  ;
\node  at    (0,3) {$w_{1}$} ;
\node  at    (0,4) {$w_{0}$}  ;
\node  at    (1,0) {$w'_{m}$}  ;
\node  at    (1,1) {$w'_{m-1}$} ;
\node  at    (1,2) {$w'_{2}$}  ;
\node  at    (1,3) {$w'_{1}$}  ;
\node  at    (1,4) {$w'_{0}$} ;
\node[fill=white] at    (0.2,0.5) {$a_{m}$}  ;
\node[fill=white] at    (0.2,2.5) {$a_{2}$}  ;
\node[fill=white] at    (0.2,3.5) {$a_{1}$}  ;
\node[fill=white] at    (1.2,0.5) {$b_{m}$}  ;
\node[fill=white] at    (1.2,2.5) {$b_{2}$}  ;
\node[fill=white] at    (1.2,3.5) {$b_{1}$}  ;
\node[fill=white] at    (0.5,0.2) {$c_{m}$}  ;
\node[fill=white] at    (0.5,2.2) {$c_{2}$}  ;
\node[fill=white] at    (0.5,3.2) {$c_{1}$}  ;
\node[fill=white] at    (0.5,4.2) {$c_{0}$}  ;
\node[fill=white] at    (0.5,1.5) {$\vdots$}  ;
\draw[->][thick,black] (0,0.2) -- (0,0.8);
\draw[->][thick,black] (0,2.2) -- (0,2.8);
\draw[->][thick,black] (0,3.2) -- (0,3.8);
\draw[->][thick,black] (1,0.2) -- (1,0.8);
\draw[->][thick,black] (1,2.2) -- (1,2.8);
\draw[->][thick,black] (1,3.2) -- (1,3.8);
\draw[->][thick,black] (0.8,0) -- (0.2,0);
\draw[->][thick,black] (0.7,1) -- (0.3,1);
\draw[->][thick,black] (0.8,2) -- (0.2,2);
\draw[->][thick,black] (0.8,3) -- (0.2,3);
\draw[->][thick,black] (0.8,4) -- (0.2,4);
\end{tikzpicture}
\end{center}
Then we have the identity $$\pi(w_{0})+\pi(w_{1})+...+\pi(w_{m})=c_{0}+\pi(w'_{1})+...+\pi(w'_{m})$$

This can be proved by induction as follows. It is clear for $m=0$, so assume its true for $m'<m$. Form the path by taking the vertex $w_{0}$ and adding the minimal incoming arrow (if we have a choice add the vertical one) and then doing the same for the initial vertex of this path. Continue until we have to choose a horizontal arrow, say we end up at $w'_{k}$. This means $a_{i} \leq c_{i-1}$ for $1 \leq i \leq k$ and $c_{k}\leq a_{k+1}$. The inequalities give $\pi(w_{i})=a_{i+1}$ for $0 \leq i \leq k-1$ and $\pi(w_{k})=c_{k}$. But the inequalities also imply $c_{i} \geq b_{i}$ for $1 \leq i \leq k$ and $c_{k+1}\leq b_{k+1}$, so $\pi(w'_{i})=b_{i}$ for $1 \leq i \leq k$ and $\pi(w'_{k+1})=c_{k+1}$. Hence $\pi(w_{0})+\pi(w_{1})+...+\pi(w_{k})=c_{0}+\pi(w'_{1})+...+\pi(w'_{k})$ as both sides are a sum over a path from $w'_{k}$ to $w_{0}$. If $k=m$ were done. If $k<m$ by the induction hypotheses we have $\pi(w_{k+1})+...+\pi(w_{m})=c_{k+1}+\pi(w'_{k+2})+...+\pi(w'_{m})$, so we get the result by adding these two equations.

Now we can use this identity to complete the proof of the lemma. By the identity we have $$\pi(v)+\sum\limits_{bel(v)}\pi(u)=c+\sum\limits_{bel(w)}\pi(u)$$ so we get  \begin{align*}\sum\limits_{bel(w)}\pi(u)-\sum\limits_{lef(w)}\pi(u)&=\sum\limits_{bel(v)}\pi(u)-c+\pi(v)-\sum\limits_{lef(w)}\pi(u)\\
&=-c+\sum\limits_{bel(v)}\pi(u)-\sum\limits_{lef(v)}\pi(u)\\
&=-c+\delta_{v}\\
&=\delta_{w}
\end{align*}
Hence horizontal induction goes through. Vertical induction is proved similarly. This completes the proof of Lemma \ref{linv}. $\square$

Now we can use Lemma \ref{linv} to compute $\langle \sum n_{ij} \alpha_{ij},\mu_{k} \rangle=\sum\limits_{l>k}n_{kl}-\sum\limits_{l<k}n_{lk}=\delta_{v_{kk}}=\lambda_{k}$. Hence $\{n_{ij}\}$ is a filling for $\lambda$.

Next we will show it is an ideal filling. We need to show that if we have a sub-diagram that looks like
\begin{center}
 \begin{tikzpicture}[thick, scale=1.3]
\node  at    (0,0) {$v$} ;
\node  at    (1,0) {$u$}  ;
\node  at    (0,1) {$w$}  ;
\node[fill=white] at    (0.5,0.13) {$b$}  ;
\node[fill=white] at    (0.13,0.5) {$a$}  ;
\draw[->][thick,black] (0.8,0) -- (0.2,0);
\draw[->][thick,black] (0,0.2) -- (0,0.8);
\end{tikzpicture}
\end{center}
then we must have $\pi(v)=\max\{\pi(u),\pi(w)\}$. First we claim $\pi(v)\geq\pi(u)$. To see this form the path $\Gamma$ by starting with $v$ and adding the minimal incoming arrow (pick the vertical one if we have a choice), and stop after the first horizontal arrow is added. Say $\Gamma$ starts at a vertex $y$. Let $\Gamma'$ be the path from $y$ to $u$ to $v$. Then all the arrows in $\Gamma$ are $\leq \pi(v)$ (by the critical conditions) and all the arrows in $\Gamma'$ are $\geq b$ (by the box-relations). Hence $\mathrm{len}(\Gamma)\pi(v) \geq \mathrm{sum}(\Gamma) =\mathrm{sum}(\Gamma') \geq \mathrm{len}(\Gamma')b\geq \mathrm{len}(\Gamma')\pi(u)$ and we are done. (Here $\mathrm{sum}(\Gamma)$ means the sum of the arrows in $\Gamma$). Similarly $\pi(v)\geq\pi(w)$.\\

Now suppose $\pi(v)>\pi(w)$ the we will show $\pi(v)=\pi(u)$. First note that $\pi(v)>\pi(w)$ forces $\pi(u)=b$ since otherwise there is a vertical arrow $d$ from $u$ to a vertex $x$ with $d<b$ and then if we let $c$ be the arrow from $x$ to $w$ then we have $c>a$ (by the box-relation) which implies $\pi(w)=a \geq \pi(v)$. Hence we have $b \geq \pi(v) \geq \pi(u) \geq b$ and we're done. Thus we do indeed get an ideal filing for $\lambda$.\\

Finally, these maps are inverse to each other because if $\{n_{ij}\}$ is an ideal filling and $\delta_{ji}=H^{h}_{ij}-H^{v}_{ij}$ then $\pi(v_{ji})=n_{ij}$ by (6.1) and (6.2), and the other way follows from Lemma \ref{linv}. Also, the claim about integrality holds since both the maps defined preserve integrality. This completes the proof of Proposition 6.2. $\square$\\

\section{A distinguished point in the Feigin Fourier Littelmann Polytope}
In this short section we show that just as the tropical critical point in the coordinates given by $\theta_{\mathcal{M}}$ gives a distinguished point in the Gelfand-Zetlin polytope, the unique ideal filling for $\lambda$ gives a distinguished point in the FFL polytope associated to $\lambda$. This is proved in Proposition 7.3 below.

The FFL polytope was introduced by Feigin, Fourier and Littelmann in [FFL]. To give the definition we first must define the Dyck paths.
\begin{defn}
A Dyck path is a sequence $(\beta(0),...,\beta(k))$ in the set of positive roots $\Delta_{G}^{+}$ such that $\beta(0)$ and $\beta(k)$ are simple roots, and if $\beta(l) = \alpha_{i\,j}$ then either $\beta(l + 1) = \alpha_{i+1\,j}$ or $\beta(l + 1) = \alpha_{i\,j+1}$.
\end{defn}
\begin{defn}
Let $\lambda=m_{1}\omega_{1}+...+m_{n}\omega_{n}$ be a dominant weight. Consider $\mathbb{R}^{ \Delta_{G}^{+}}$ and let $n_{\beta}$ be the coordinate corresponding to $\beta\in\Delta_{G}^{+}$. Then the FFL polytope for $\lambda$ is defined inside $\mathbb{R}^{ \Delta_{G}^{+}}$ by the following inequalities: $$n_{\beta}\geq 0 \quad\text{for all}\quad \beta\in \Delta_{G}^{+}$$ and $$n_{\beta(0)}+...+n_{\beta(k)} \leq m_{i}+...+m_{j}$$ for all Dyck paths $(\beta(0),..,\beta(k))$ such that $\beta(0)=\alpha_{i}$ and $\beta(k)=\alpha_{j}$.
\end{defn}
\begin{prop}
Let $(n_{ij})_{i<j}$ be an ideal filling for $\lambda$. Then the point in $\mathbb{R}^{\Delta_{G}^{+}}$, defined by setting $n_{\beta}=n_{ij}$ for $\beta=\alpha_{ij}$, lies in the FFL polytope for $\lambda$.
\end{prop}
\noindent\textit{Proof:} We must show that the point satisfies the inequalities defining the FFL polytope. We will need the following formula for the pairing of a string of roots with a string of coroots $$\langle \alpha_{i}+..+\alpha_{j},\alpha_{k}^{\vee}+..+\alpha_{l}^{\vee} \rangle=\bm{\delta_{i\,k}}+\bm{\delta_{j+1\,k+1}}-\bm{\delta_{i\,k+1}}-\bm{\delta_{j+1\,k}}$$ where $\bm{\delta_{p\,q}}=\begin{cases} 1 \quad \text{if}\quad p=q\\ 0 \quad \text{if}\quad p\neq q \end{cases}$ (i.e. the Kronecker delta). This follows from the formula for the pairing of a root with a coroot $$\langle \alpha_{p},\alpha_{q}^{\vee} \rangle= \begin{cases} 2 \quad\quad \text{if}\quad p=q\\ -1 \quad \text{if}\quad \vert p- q\vert=1 \\0\quad \quad \text{otherwise}\end{cases}$$.

First we show $n_{ij}\geq 0$ for all $i<j$. By the ideal property we just need to show $n_{ii+1}\geq 0$ for $1\leq i \leq n-1$. Let $n_{kk+1}$ be the minimal element of $\{ n_{ii+1}:1\leq i \leq n-1\}$. Then
 \begin{align*}
 m_{k}&=\langle \lambda , \alpha_{k}^{\vee} \rangle \\&=\sum\limits_{i \leq j} n_{i\,j+1}\langle \alpha_{i}+..+\alpha_{j},\alpha_{k}^{\vee} \rangle \\&=\sum\limits_{k\leq j}n_{k\,j+1}+\sum\limits_{i\leq k}n_{i\,k+1}-\sum\limits_{k+1\leq j}n_{k+1\,j+1}-\sum\limits_{i\leq k-1}n_{i\,k}\\
 &=2n_{k\,k+1}+ \sum\limits_{k+1\leq j}n_{k\,j+1}+\sum\limits_{i\leq k-1}n_{i\,k+1}-\sum\limits_{k+1\leq j}n_{k+1\,j+1}-\sum\limits_{i\leq k-1}n_{i\,k} \\
 &=2n_{k\,k+1}+ \sum\limits_{k+1\leq j}(n_{k\,j+1}-n_{k+1\,j+1})+\sum\limits_{i\leq k-1}(n_{i\,k+1}-n_{i\,k})\\
 &=2n_{k\,k+1}\end{align*}since $n_{k\,k+1}$ being minimal implies $n_{i\,k}=\max\limits_{i\leq p\leq k-1}\{ n_{p\,p+1} \}=\max\limits_{i\leq p\leq k}\{ n_{p\,p+1} \}=n_{i\,k+1}$. Hence $\lambda$ being dominant gives $n_{k\,k+1}=m_{k}/2 \geq0$. Thus the first set of inequalities hold.\\

Now let $(\beta(0),..,\beta(t))$ be a Dyck path with $\beta(0)=\alpha_{k}$ and $\beta(t)=\alpha_{l}$. We need to show $$n_{\beta(0)}+...+n_{\beta(k)} \leq m_{i}+...+m_{j}.$$ Similar to above, we have \begin{align*} m_{k}+...+m_{l}&=\langle \lambda , \alpha_{k}^{\vee}+...+ \alpha_{l}^{\vee} \rangle \\
 &=\sum\limits_{i \leq j} n_{i\,j+1}\langle \alpha_{i}+..+\alpha_{j},\alpha_{k}^{\vee}+...+ \alpha_{l}^{\vee}  \rangle \\
 &=\sum\limits_{k\leq j}n_{k\,j+1}+\sum\limits_{i\leq l}n_{i\,l+1}-\sum\limits_{l+1\leq j}n_{l+1\,j+1}-\sum\limits_{i\leq k-1}n_{i\,k}\\
 &=\sum\limits_{k\leq j\leq l}n_{k\,j+1}+\sum\limits_{l+1\leq j}(n_{k\,j+1}-n_{l+1\,j+1})+\sum\limits_{k\leq i\leq l}n_{i\,l+1}+\sum\limits_{i\leq k-1}(n_{i\,l+1}-n_{i\,k}) \\
 &\geq \sum\limits_{k\leq j\leq l}n_{k\,j+1}+\sum\limits_{k\leq i\leq l}n_{i\,l+1}\\
 &=n_{k\,l+1}+\sum\limits_{\beta \in \mathcal{P}}n_{\beta}\\
 &\geq \sum\limits_{\beta \in \mathcal{P}}n_{\beta}
  \end{align*}
where $\mathcal{P}$ is the Dyck path given by $\alpha_{k\,k+1},\alpha_{k\,k+2},...,\alpha_{k\,l+1},\alpha_{k+1\,l+1},...,\alpha_{l\,l+1}$. Now the ideal property implies that for any Dyck path from $\alpha_{k}$ to $\alpha_{l}$, say $\mathcal{P}'$, we have that $\sum\limits_{\beta \in \mathcal{P}'}n_{\beta}\leq \sum\limits_{\beta \in \mathcal{P}}n_{\beta}$. Hence the second set of inequalities hold, and the point lies in the FFL polytope. $\square$\\

\section{A canonical section in $H^{0}(G/B,\mathcal{L}_{\lambda})$ for $\lambda\in \mathcal{P}$}
In this section we define a canonical section in $H^{0}(G/B,\mathcal{L}_{\lambda})$ for certain $\lambda$. First consider the case when $\lambda=2 \rho$. In this case $\mathcal{L}_{2 \rho}$ is the anti-canonical bundle of $G/B$. Now, there exists a special non-vanishing meromorphic top-form $\omega$ on $G/B$. It is the unique (up to scalar) meromorphic differential form on $G/B$ with simple poles exactly along the divisor given by the union of all the Schubert divisors and all the opposite Schubert divisors (see [Lam]). This form was first introduced in [R] where it was defined as a natural generalisation of the unique torus-invariant volume form on a torus inside a toric variety, and in fact this definition pins it down uniquely up to sign. Now if we take the inverse of $\omega$, we get a special global section of the anti-canonical bundle of $G/B$, and thus a distinguished vector in the representation $H^{0}(G/B,\mathcal{L}_{2\rho})$. We will generalise this to find a special section in $H^{0}(G/B,\mathcal{L}_{\lambda})$ for certain $\lambda$.\\

Let $\bigwedge^{N_{P}}\mathcal{T}^{*}_{G/P}$ be the canonical line bundle on $G/P$ and $\bigwedge^{N_{P}}\mathcal{T}_{G/P}$ the anti-canonical line bundle. There is a natural choice of holomorphic volume form (defined up to sign) on $\mathcal{R}_{w_{P},w_{0}}$ defined in [R,§7]. We denote it by $\omega_{P}$ and it is given by $$\omega_{P}=\frac{du_{k_{1}}}{u_{k_{1}}}\wedge ...\wedge \frac{du_{k_{N_{P}}}}{u_{k_{N_{P}}}}$$ in any of the toric charts $\textbf{y}^{P}_{\textbf{i}}$ from section 3. This volume form extends to give a non-vanishing meromorphic top form on $G/P$ which we also denote by $\omega_{P}$. Let $\omega_{P}^{-1} \in \bigwedge^{N_{P}}\mathcal{T}_{G/P} $ be the dual section to $\omega_{P}\in \bigwedge^{N_{P}}\mathcal{T}^{*}_{G/P}$. Then $\omega_{P}^{-1}$ is a global section of $\bigwedge^{N_{P}}\mathcal{T}_{G/P}$.\\

Next consider the map $\pi_{P}:G/B \rightarrow G/P$ and look at the pull back of $\bigwedge^{N_{P}}\mathcal{T}_{G/P}$ under $\pi_{P}$. The fiber at the basepoint is isomorphic to $\mathfrak{g}/\mathfrak{p}$, so the $T$-weights are $\Delta_{-}^{G}\backslash \Delta_{-}^{L}$, where $L$ is the Levi of $P$. Then the sum of these $T$-weights is equal to $-\lambda_{P}$, since we have $\lambda_{P}=\sum\limits_{\Delta_{+}^{G}\backslash \Delta_{+}^{L}} \alpha$. Hence the pullback of the anti-canonical bundle is the line bundle $\mathcal{L}_{\lambda_{P}}$.\footnote{This is nicely explained in [Kn]}\\

Now for $\lambda\in \mathcal{P}$ we can define a section analogous to the anti-canonical section given above.

\begin{defn}(\textit{The section $\omega_{\lambda}^{-1}$}) Given $\lambda\in\mathcal{P}$ there is a unique way to write $\lambda=\sum\limits_{P\subset G}c_{P}\lambda_{P}\text{ where }c_{P}\in \mathbb{Z}_{\geq0}$ and the set of parabolics $\{P:c_{P}\neq 0\}$ form a chain w.r.t inclusion. Now we have a surjection $$\bigotimes\limits_{P} H^{0}(G/B,\mathcal{L}_{\lambda_{P}})^{\otimes c_{P}} \rightarrow H^{0}(G/B,\mathcal{L}_{\lambda})$$ and we define $\omega_{\lambda}^{-1}$ to be the image of $\bigotimes\limits_{P}(\pi_{P}^{*} \omega_{P}^{-1})^{\otimes c_{P}}$ under this map.\\
\end{defn}

\begin{remark}
Here are the dominant integral weights for $SL_3$ with those lying in $\mathcal{P}$ circled.\\
\end{remark}
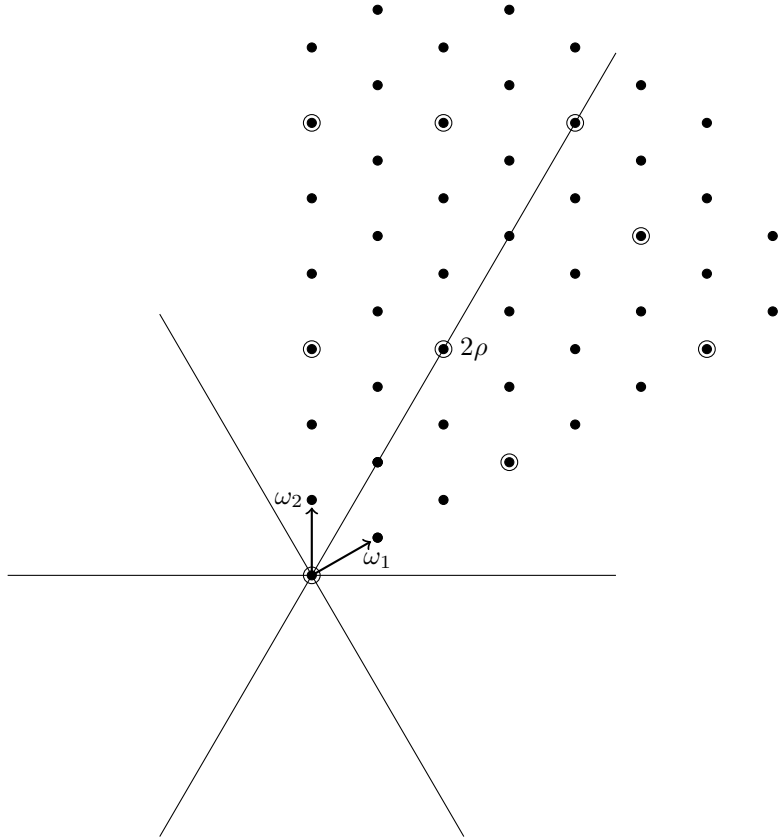
\begin{center}
 \begin{tikzpicture}[scale=1]
\draw [fill] (0,1)  circle [radius=0.06];
\draw [fill]  (0.866,0.5)  circle [radius=0.06];
\draw [fill]  (0,1)++(0.866,0.5)  circle [radius=0.06];
\draw [fill]  (0,2*1)  circle [radius=0.06];
\draw [fill]  (0,1)++(2*0.866,2*0.5)  circle [radius=0.06];
\draw [fill]  (0,3)  circle [radius=0.06];
\draw [fill]  (0,4) circle [radius=0.06];
\draw [fill]  (0,5*1)  circle [radius=0.06];
\draw [fill]  (0,6*1)  circle [radius=0.06];
\draw [fill]  (0,7)  circle [radius=0.06];
\draw [fill]  (0,1)++(0.866,0.5)  circle [radius=0.06];
\draw [fill]  (0,2*1)++(0.866,0.5)  circle [radius=0.06];
\draw [fill]  (0,3)++(0.866,0.5)  circle [radius=0.06];
\draw [fill]  (0,4)++(0.866,0.5)  circle [radius=0.06];
\draw [fill]  (0,5*1)++(0.866,0.5)  circle [radius=0.06];
\draw [fill]  (0,6*1)++(0.866,0.5)  circle [radius=0.06];
\draw [fill]  (0,1)++(0.866,0.5)  circle [radius=0.06];
\draw [fill]  (0,2)++(2*0.866,2*0.5)  circle [radius=0.06];
\draw [fill]  (0,3)++(2*0.866,2*0.5)  circle [radius=0.06];
\draw [fill]  (0,4)++(2*0.866,2*0.5)  circle [radius=0.06];
\draw [fill]  (0,5*1)++(2*0.866,2*0.5)  circle [radius=0.06];
\draw [fill]  (0,1)++(3*0.866,3*0.5)  circle [radius=0.06];
\draw [fill]  (0,2)++(3*0.866,3*0.5)  circle [radius=0.06];
\draw [fill]  (0,3)++(3*0.866,3*0.5)  circle [radius=0.06];
\draw [fill]  (0,4)++(3*0.866,3*0.5)  circle [radius=0.06];
\draw [fill]  (0,1)++(4*0.866,4*0.5)  circle [radius=0.06];
\draw [fill]  (0,2)++(4*0.866,4*0.5)  circle [radius=0.06];
\draw [fill]  (0,3)++(4*0.866,4*0.5)  circle [radius=0.06];
\draw [fill]  (0,1)++(5*0.866,5*0.5)  circle [radius=0.06];
\draw [fill]  (0,2)++(5*0.866,5*0.5)  circle [radius=0.06];
\draw [fill]  (0,1)++(6*0.866,6*0.5)  circle [radius=0.06];
\draw [fill]  (0.866,0.5)  circle [radius=0.06];
\draw [fill]  (2*0.866,2*0.5)  circle [radius=0.06];
\draw [fill]  (3*0.866,3*0.5)  circle [radius=0.06];
\draw [fill]  (4*0.866,4*0.5)  circle [radius=0.06];
\draw [fill]  (5*0.866,5*0.5)  circle [radius=0.06];
\draw [fill]  (6*0.866,6*0.5)  circle [radius=0.06];
\draw [fill]  (7*0.866,7*0.5)  circle [radius=0.06];

\draw [fill]  (0,6)++(2*0.866,2*0.5)  circle [radius=0.06];
\draw [fill]  (0,5)++(3*0.866,3*0.5)  circle [radius=0.06];
\draw [fill]  (0,4)++(4*0.866,4*0.5)  circle [radius=0.06];
\draw [fill]  (0,3)++(5*0.866,5*0.5)  circle [radius=0.06];
\draw [fill]  (0,2)++(6*0.866,6*0.5)  circle [radius=0.06];
\draw [fill]  (0,5)++(4*0.866,4*0.5)  circle [radius=0.06];
\draw [fill]  (0,4)++(5*0.866,5*0.5)  circle [radius=0.06];

\draw [fill]  (0,7)++(1*0.866,1*0.5)  circle [radius=0.06];
\draw [fill]  (0,6)++(3*0.866,3*0.5)  circle [radius=0.06];
\draw [fill]  (0,3)++(6*0.866,6*0.5)  circle [radius=0.06];
\draw [fill]  (0,1)++(7*0.866,7*0.5)  circle [radius=0.06];

\draw [circle]  (0,0)  circle [radius=0.11];
\draw [fill]  (0,0.0)  circle [radius=0.06];

\draw [circle]  (0,3)  circle [radius=0.11];
\draw [circle]  (0,6)  circle [radius=0.11];
\draw [circle]  (3*0.866,3*0.5)   circle [radius=0.11];
\draw [circle]  (6*0.866,6*0.5)   circle [radius=0.11];
\draw [circle]  (0,2)++(2*0.866,2*0.5)   circle [radius=0.11];
\draw [circle]  (0,2)++(2*0.866,2*0.5) ++ (0,3) circle [radius=0.11];
\draw [circle]  (0,2)++(2*0.866,2*0.5) ++ (3*0.866,3*0.5)  circle [radius=0.11];
\draw [circle]  (0,4)++(4*0.866,4*0.5)   circle [radius=0.11];

\draw[->][thick,black] (0,0) -- (0,0.9*1);
\draw[->][thick,black] (0,0) -- (0.9*0.866,0.9*0.5);
\node  at    (-0.3,1) {$\omega_{2}$} ;
\node  at    (0.866,0.2) {$\omega_{1}$} ;

\node  at   (2*1.066,2*1.5)  {$2\rho$} ;

\draw[-][thin,black] (-4,0) -- (4,0);
\draw[-][thin,black] (-4*0.5,-4*0.866) -- (8*0.5,8*0.866);
\draw[-][thin,black] (-4*0.5,4*0.866) -- (4*0.5,-4*0.866);

\end{tikzpicture}
\end{center}
\bigbreak

\section{Conjecture about $\omega_{\lambda}^{-1}$ and $p_{\lambda}$}

Recall that given $\textbf{i}$ a reduced expression of $w_0$, we have two maps
\begin{align*}
\nu^{\vee}_{\lambda,\textbf{i}}:Z_{t^\lambda}(\mathcal{K}_{>0})^+\longrightarrow \mathbb{R}^N\\
\nu_{\lambda,\textbf{i}}:H^{0}(G/B,\mathcal{L}_{\lambda})\longrightarrow \mathbb{R}^N
\end{align*}
both of whos image lies in the string polytope $\mathrm{String}_{\bf{i}}(\lambda)$.
For $\lambda\in\mathcal{P}$ we have constructed a section $\omega_{\lambda}^{-1}\in H^{0}(G/B,\mathcal{L}_{\lambda})$, and for $\lambda\in P^+$ we have constructed a point $p_{\lambda}\in Z_{t^\lambda}(\mathcal{K}_{>0})^+$, which is integral iff $\lambda\in\mathcal{P}$. Then we conjecture that
\begin{conj}
Given $\lambda\in\mathcal{P}$, we have $\nu_{\lambda,\textbf{i}}(\omega_{\lambda}^{-1})=\nu^{\vee}_{\lambda,\textbf{i}}(p_{\lambda})$ for all $\textbf{i}$ reduced expressions of $w_0$.
\end{conj}

\begin{example}
Let $n=3$. Let $\lambda=2 \omega_{1}+5 \omega_{2}$. Note $\lambda=2 \rho+3 \omega_{2}=\lambda_{B}+\lambda_{P}$ where $P=\langle B,\dot{s}_{1} \rangle$. Let's take $\textbf{i}=(212)$. First we will compute $\nu_{\lambda,\textbf{i}}(\omega_{\lambda}^{-1})$. In the coordinates on $G/B$ given by $\textbf{y}_{\textbf{i}}^{B}$, we have $$\omega_{B}^{-1}=u_{1} u_{2} u_{3} \partial_{u_{1}} \wedge \partial_{u_{2}} \wedge \partial_{u_{3}}$$ where $\partial_{u}$ is shorthand for $\frac{\partial}{\partial u}$. Now we claim that a lowest weight section is given by the pullback of $\frac{1}{u_{2}} \partial_{u_{1}} \wedge \partial_{u_{2}} \wedge \partial_{u_{3}}$. For this we need to show its $T$-weight is $-\lambda_{B}$ which is the lowest weight of $V_{\lambda_{B}}^{*}$, i.e. $$\begin{pmatrix}t_{1}&&\\&t_{2}&\\&&t_{3}\end{pmatrix}\cdot\frac{1}{u_{2}} \partial_{u_{1}} \wedge \partial_{u_{2}} \wedge \partial_{u_{3}}=(-\lambda_{B})(t)\cdot \frac{1}{u_{2}} \partial_{u_{1}} \wedge \partial_{u_{2}} \wedge \partial_{u_{3}}$$ We compute the $T$-weight by solving $$ \begin{pmatrix}t_{1}&&\\&t_{2}&\\&&t_{3}\end{pmatrix}^{-1} \textbf{y}_{2}(u_{1})  \textbf{y}_{1}(u_{2}) \textbf{y}_{2}(u_{3})B =  \textbf{y}_{2}(v_{1}) \textbf{y}_{1}(v_{2}) \textbf{y}_{2}(v_{3})B$$ to get $u_{1}=\frac{t_{3}v_{1}}{t_{2}}$ and $u_{2}=\frac{t_{2}v_{2}}{t_{1}}$ and $u_{3}=\frac{t_{3}v_{3}}{t_{2}}$ whence $\frac{t_{2}}{t_{1}}\det ( \frac{\partial(u_{1},u_{2},u_{3})}{\partial(v_{1},v_{2},v_{3})})=\frac{t_{3}^{2}}{t_{1}^{2}}=-\lambda_{B}(t)$. So the $T$-weight is as claimed. Let $$f_{B}:=\frac{\pi_{B}^{*}\omega_{B}^{-1}}{\pi_{B}^{*}\frac{1}{u_{2}} \partial_{u_{1}} \wedge \partial_{u_{2}} \wedge \partial_{u_{3}} }\in\mathbb{C}[\mathcal{R}_{e,w_{0}}].$$ Then $f_{B}$ is given by the pullback of $u_{1}u_{2}^{2}u_{3}$. In the coordinates on $G/B$ given by $\textbf{y}_{\textbf{i}}$ this pullback is simply given by $x_{1}x_{2}^{2}x_{3}$.

Similarly, in the coordinates on $G/P$ given by $\textbf{y}_{\textbf{i}}^{P}$, we have $\omega_{P}^{-1}=u_{1} u_{3} \partial_{u_{1}} \wedge \partial_{u_{3}} $. A lowest weight section is given by the pullback of $\partial_{u_{1}} \wedge \partial_{u_{3}}$. For this we need to show its $T$-weight is $-\lambda_{P}$ which is the lowest weight of $V_{\lambda_{P}}^{*}$, i.e.  $$\begin{pmatrix}t_{1}&&\\&t_{2}&\\&&t_{3}\end{pmatrix}\cdot\partial_{u_{1}} \wedge \partial_{u_{3}}=(-\lambda_{P})(t)\cdot \partial_{u_{1}} \wedge \partial_{u_{3}}$$ In this case, computing the $T$-weight is done by solving  $$ \begin{pmatrix}t_{1}&&\\&t_{2}&\\&&t_{3}\end{pmatrix}^{-1} \textbf{y}_{2}(u_{1}) \dot{s}_{1} \textbf{y}_{2}(u_{3})P =  \textbf{y}_{2}(v_{1})\dot{s}_{1} \textbf{y}_{2}(v_{3})P$$ to get $u_{1}=\frac{t_{3}v_{1}}{t_{2}}$ and $u_{3}=\frac{t_{3}v_{3}}{t_{1}}$ whence $\det ( \frac{\partial(u_{1},u_{3})}{\partial(v_{1},v_{3})})=\frac{t_{3}^{2}}{t_{1}t_{2}}=-\lambda_{P}(t)$. Now let $$f_{P}:=\frac{\pi_{P}^{*}\omega_{P}^{-1}}{\pi_{P}^{*}\partial_{u_{1}} \wedge \partial_{u_{3}}}\in\mathbb{C}[\mathcal{R}_{e,w_{0}}].$$ Then $f_{P}$ is given by the pullback of $u_{1}u_{3}$. With respect to the chart $\textbf{y}_{\textbf{i}}$ on $G/B$ and $\textbf{y}_{\textbf{i}}^{P}$ on $G/P$, the map $\pi_{P}$ is given by $$\pi_{P}:(x_{1},x_{2},x_{3}) \mapsto (x_{1}+x_{3},-x_{2}x_{3}).$$ Hence, in the coordinates on $G/B$ given by $\textbf{y}_{\textbf{i}}$ we have $f_{P}=-x_{2}x_{3}^{2}-x_{1}x_{2}x_{3}$.

Since the lowest weight section of a tensor product of representations is just the tensor product of the lowest weight sections, the image of $\omega_{\lambda}^{-1}$ in $\mathbb{C}[\mathcal{R}_{e,w_{0}}]$ is just gotten by multiplying $f_{B}$ and $f_{P}$. Hence we get $-x_{1}x_{2}^{3}x_{3}^{3}-x_{1}^{2}x_{2}^{3}x_{3}^{2}$. The lexicographically maximal term of this is $-x_{1}^{2}x_{2}^{3}x_{3}^{2}$, so we get the point $\nu_{\lambda,\textbf{i}}(\omega_{\lambda}^{-1})=(2,3,2)\in \mathbb{Z}^{3}$.\\

Now we compute $\nu^{\vee}_{\lambda,\textbf{i}}(p_{\lambda})$. It is easiest to compute in the special chart and then change coordinates. The transition map $(\tilde{\textbf{x}}_{-\textbf{i}}^{\vee} )^{-1}\circ \theta_{\mathcal{M}} $ is given by $$(a,b,c,d,e,f)\mapsto \left( \begin{pmatrix}acdf&&\\&df&\\&&1\end{pmatrix} ,(\frac{ad}{b+d},ab,b+d)\right)$$ Also note that to compute $\nu^{\vee}_{\lambda,\textbf{i}}(p_{\lambda})$ we just need to compute the valuation of the critical point in the special chart and then apply the tropicalisation of the transition map above. In the coordinates given defined by the chart $\theta_{\mathcal{M}}$ it is easy to compute the tropical critical point. It is given by the unique way to fill in the diagram with $\tilde{\lambda}=(7,5,0)$ on the diagonal such that the tropical critical point relations are satisfied.
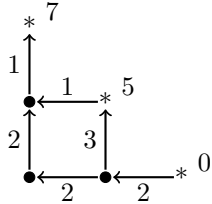
\begin{center}
 \begin{tikzpicture}
\draw [fill] (0.5,0.5) circle [radius=0.07];
\draw [fill] (0.5,1.5) circle [radius=0.07];
\draw [fill] (1.5,0.5) circle [radius=0.07];
\node  at    (0.5,2.5) {*} ;
\node  at    (1.5,1.5) {*}  ;
\node  at    (2.5,0.5) {*}  ;

\node[fill=white] at    (0.3,1) {2}  ;
\node[fill=white] at    (0.3,2) {1}  ;
\node[fill=white] at    (1.3,1) {3}  ;
\node[fill=white] at    (1,0.3) {2}  ;
\node[fill=white] at    (1,1.7) {1}  ;
\node[fill=white] at    (2,0.3) {2}  ;

\node[fill=white] at    (0.8,2.7) {7}  ;
\node[fill=white] at    (1.8,1.7) {5}  ;
\node[fill=white] at    (2.8,0.7) {0}  ;

\draw[->][thick,black] (0.5,0.6) -- (0.5,1.4);
\draw[->][thick,black] (0.5,1.6) -- (0.5,2.4);
\draw[->][thick,black] (1.5,0.6) -- (1.5,1.4);
\draw[->][thick,black] (2.4,0.5) -- (1.6,0.5);
\draw[->][thick,black] (1.4,0.5) -- (0.6,0.5);
\draw[->][thick,black] (1.4,1.5) -- (0.6,1.5);
\end{tikzpicture}
\end{center}
Transforming this point to the chart $ \tilde{\textbf{x}}_{-\textbf{i}}^{\vee} $ using the tropicalisation of the transition map above gives the point $(1+3-\min(2,3),1+2,\min(2,3))=(2,3,2)\in \mathbb{Z}^{3}$ which agrees with above.

\end{example}

\end{document}